\documentclass[draft]{amsart}

\usepackage{amssymb}
\usepackage{url}
\usepackage{amscd}

\numberwithin{equation}{section}
\theoremstyle{plain}
\newtheorem{theor}{Theorem}

\newtheorem{Th}[equation]{Theorem}
\newtheorem{Prop}[equation]{Proposition}
\newtheorem{Le}[equation]{Lemma}
\newtheorem{Cor}[equation]{Corollary}

\theoremstyle{definition}
\newtheorem{definition}[equation]{Definition}
\newtheorem{example}[equation]{Example}
\newtheorem{remark}[equation]{Remark}

\newcommand{\inter}{\operatorname{Int}}
\newcommand{\sha}{\operatorname{Shadow}}

\newcommand{\Z}{\mathbb Z}
\newcommand{\R}{\mathbb R}
\newcommand{\Q}{\mathbb Q}
\newcommand{\N}{\mathbb N}

\begin{document}

\date{Received August 11, 2003; received in final form January 11, 2004} 

\title[Isometries of $CAT(0)$-spaces as maps preserving diagonal tube]{A characterization of isometries of $CAT(0)$-space as maps preserving diagonal tube}

\author{Pavel D. Andreev}

\address{Pomor State University, Arkhangelsk, Russia}

\email{andreev@math.pomorsu.ru}


\subjclass[2000]{Primary 57C23, 57C70}

{\sloppy
\begin{abstract}{We give positive answers for  questions by Berestovski\v\i. Namely, we prove that every bijection of locally compact geodesically complete and connected at infinity $CAT(0)$-space $X\/$ onto itself preserving some fixed distance or satellite relations is an isometry of this space.
The proof of this theorem is based on another result stated by Berestovski\v\i~ as a problem: the metric of the space $X\/$ may be recovered from its diagonal tube corresponding to an arbitrary number $r \, >\, 0$.}
\end{abstract}

}

\maketitle

\vspace{-5mm}



      

\section*{Introduction}

{\sloppy
V.N.Berestovski\u\i~  in \cite{Berestovskii} has established following characterisation for isometries of Aleksandrov spaces of curvature negatively bounded above:

}
\begin{theor}[Berestovski\u\i]\label{isom}
Let $X\/$ be a locally compact geodesically complete $CAT(\kappa)$-space, $\kappa \, <\, 0$, in which all spheres are path connected. Then every bijection $f\/$ of $X\/$ onto itself such that both $f\/$ and $f^{-1}\/$ map any closed ball of some fixed radius $r\, >\,0$ onto some ball of radius $r\/$ is isometry of $X$.
\end{theor}

The key point of the proof is

\begin{theor}[Berestovski\u\i]
Let $X\/$ be as in Theorem \ref{isom}  and  $V \subset X \times X\/$ be its diagonal tube corresponding to a number $r\, >\,0$. Then the metric of $X\/$ is uniquely determined by $V$.
\end{theor}

A \emph{diagonal tube\/} $V\/$ of a metric space $X\/$ corresponding to $r\, >\,0$ is by definition a set
\[V := \{ (x,y) \in X \times X |\quad | xy| \, \le\, r\} \subset X \times X\, ,\]
where $|xy|\/$ is a distance between points $x,y \in X$. For a set $V\/$ we put 
\[\partial V := \{ (x,y) \in X \times X |\quad | xy|\, =\, r\} \subset X \times X\, ,\]
and
\[{\rm Int}V := \{ (x,y) \in X \times X |\quad | xy| \, <\, r\} \subset X \times X\, .\]

A question whether analogous theorems are true in the case $\kappa\, =\, 0\/$ is raised in \cite{Berestovskii}. We give an affirmative answer to this question. Namely we prove

\begin{Th}\label{main}
Let $(X, d)\/$ be a locally compact geodesically complete and connected at infinity $CAT(0)$-space, $f\,\colon X \, \to\, X$ --- bijection and  $V \subset X \times X\/$ be diagonal tube of space $X\/$ corresponding to a number $r\, >\,0$. Set a map $\phi := f\times f \,\colon X\times X \, \to\, X \times X\/$ by $\phi(x_1, x_2) := (f(x_1), f(x_2))$. Then following statements are equivalent:
\begin{enumerate}
\item $\phi(V)\, =\, V$,
\item $\phi(\partial V)\, =\, \partial V$,
\item $\phi(\inter V)\, =\, \inter V$,
\item $f\/$ is an isometry of $X\/$ onto itself.
\end{enumerate}
\end{Th}

\begin{remark}
When we say that a map $\phi\/$ moves $V\/$ onto $V\/$ this means that for every pair $x,y \in X\/$ with $d(x,y) \, \le\, r\/$ inequality $d(\phi(x), \phi(y)) \, \, \le\, r\/$ holds, but not necessarily $d(\phi(x), \phi(y))\, =\, d(x, y)$. The last equality is the target of the theorem. 
\end{remark}

\begin{remark}
Connectedness at infinity of space $X\/$ means that the complement $X\setminus B\/$ of every metric ball $B \subset X\/$ is path connected. It was shown in \cite{Berestovskii} that connectedness at infinity is equivalent to the condition of path connectedness for every sphere in $X$.
\end{remark}

Following \cite{Berestovskii} we base the proof on the uniqueness of the metric on $X\/$ with prescribed properties and given diagonal tube. We formulate this statement as following.

\begin{Th}\label{tube}
Let the space $(X, d)\/$ be as in theorem \ref{main} and the number $r\, >\,0\/$ is fixed. Then every metric $d'\/$ on $X\/$ such that $(X, d')\/$ is locally compact, geodesically complete $CAT(0)$-space, coincides with $d$:  for any $x,y\in X\, d'(x,y)=d(x,y)$, iff any one of following equivalent conditions hold:
\begin{enumerate}
\item $\forall x,y\in X\quad d(x,y) \, \le\, r \, \Leftrightarrow\, d'(x,y)\, \le\, r$;
\item $\forall x,y\in X\quad d(x,y)\, =\, r\, \Leftrightarrow\, d'(x,y)\, =\, r$;
\item $\forall x,y\in X\quad d(x,y) \, <\, r\, \Leftrightarrow\, d'(x,y) \, <\, r$.
\end{enumerate}
\end{Th}

Theorem \ref{main} is direct consequence of theorem \ref{tube}: it sufficies to take the metric $d'\, =\, f^\ast d$:
\[d'(x,y)\, =\, d(f(x), f(y))\, ,\]
and apply equality $f^\ast d (x,y) := d(x,y)\/$ as claim of theorem \ref{tube}.
 
For simplicity, from now on we will assume $r=1$, but conveniently, the notation $r\/$ will be kept in some notions, such as \emph{$r$-sequences\/} etc.

The author have considered a case of Busemann spaces (\cite{Busemann}), i.e.  locally compact complete metric spaces with intrinsic metric and properties of nonbranching  and local extendability of minimizing segments in earlier paper (\cite{andreev}). It was shown in  \cite{berestov-2} that every $CAT(0)$-space which is a Busemann space is really a Riemannian manifold with continuous metric tensor relatively distance coordinates. Thus the connectedness at infinity of a space $X\/$ is equivalent to  the estimation $n\, >\,1\/$ for its topological dimension $n :=  \, \operatorname{TopDim}(X)\/$  in this case. The proof of the theorem \ref{tube} in the present paper is based on methods developed in \cite{andreev} with necessary modification.

Problems resolved in the paper are conjugated with the question suggested by A.D.Alek\-sand\-rov in 1960-s:
\begin{itemize}
\item 
Under what conditions is a map of a metric space into itself preserving a fixed distance (for example, distance 1) an isometry of this metric space?
\end{itemize}

A number of similar questions are partially resolved by Aleksandrov himself. For example, several theorems presenting sufficient conditions for map of classical spaces of constant curvature preserving congruence to be isometry are proved in \cite{A.D.1973}. 
Beckman and Quarles have proved the version of theorem \ref{main} for maps of Euclidean spaces in  \cite{Beck-Quart}, Kuzminykh --- for Lobachevskii space in \cite{Kuzm}.

The paper consists of three sections. The first section is introductory and contains a list of basic notions and facts used by author as preliminaries and the concept of parallel-equivalence for $r$-sequences. We introduse two types of $r$-sequences according to behavior of corresponding classes of parallel-equivalence. In the rest sections we show that the metric of any geodesic in $X\/$ can be recovered from $V$. In the Section 2 we study geodesics bounding flat strip, and in Section 3 geodesics of strictly rank one. Combining two mentioned situations we get whole proof of Theorem \ref{tube}.

{\bf Acknowledgment} I'm very grateful to V.N.Berestovski\v\i~for support and useful consultations given during my work on the paper. I would like to thank members of geometrical seminar of Steklov PDMI RAS (St. Petersburg) for attention and number of important remarks. 
Especially I want to thank S.V.Buyalo, for suggesting an idea simplifying the proof of existence theorem  for scissors.

\section{Preliminaries }

\subsection{$CAT(\kappa)$-spaces} \label {catk}

Main definitions and properties of Aleksandrov spaces with curvature bounded above and so called $CAT(\kappa)$-spaces may be found in \cite{B}, \cite{BH} or \cite{Bu}.

Let $(X,d)\/$ be a metric space. The distance $d\/$ between points $x,y \in X\/$ conveniently will be denoted as $|xy| := d(x,y)$.
An open ball of radius $\rho\/$ centered at point  $x\in X\/$ is denoted as  $B(x, \rho)$, corresponding closed ball --- as  $\overline{B}(x, \rho)$, the boundary sphere --- as $S(x, \rho)$.
For any subset $A\subset X\/$ and any $\epsilon\, >\,0$, the $\epsilon$-neighborhood of $A\/$ is  
\[N_\epsilon(A):=\{x\in X|\, d(x, a) \, <\, \epsilon\quad \text{for some } a\in A\}\, .\]
For any two closed subsets  $A, B\subset X$,  the \emph{Hausdorff distance\/} between $A\/$  and  $B\/$ is 
\[d_H(A, B):=\inf \{\,\epsilon\,|\, A\subset N_\epsilon(B),\, B\subset N_\epsilon(A)\};\]
$d_H(A, B)\/$  is defined to be $\infty\/$ if there is no   $\epsilon\, >\,0\/$   with  $A\subset N_\epsilon(B)\/$  nor $B\subset N_\epsilon(A)$.

A \emph{geodesic\/}  in $X\/$   is a continuous map $c\,\colon I\, \rightarrow\, X\/$  from an interval  $I\subset \R\/$ into $X$,  such that for any point $t\in I\/$ there exists a neighborhood $U\/$ of $t\/$ with  $d(c(s_1),c(s_2))=|s_1-s_2|\/$ for   all  $s_1,s_2 \in U$. If one can take $U\, =\, I$, then such geodesic is said  to be \emph{minimizing}.  The image of a geodesic or a minimizing will also be called a geodesic or a minimizing.   When $I\/$ is a closed interval $[a,b]\subset \R$, we say  that $c\/$ is a \emph{geodesic segment\/} of length $b-a\/$ and  $c\/$ \emph{connects\/}  $c(a)\/$  and  $c(b)$.  If $I\, =\, \R$, we say that $c\/$ is a complete geodesic. A metric space $X\/$ is called a \emph{geodesic metric space\/} if for any  two points  $x,y\in X\/$ there is a    minimal geodesic segment connecting them.  A geodesic metric space is called \emph{geodesically complete\/} if every geodesic segment is contained in some complete geodesic (not necessarily unique).

It follows from Hopf-Rinow  theorem (cf. \cite[Ch.1, Theorem 2.3]{Bu}) that geodesically complete locally compact space $X\/$ is \emph{proper\/} or \emph{finitely compact}, i.e. every its closed ball is compact.

For $\kappa \in \R\/$ we let $M_\kappa\/$ be the \emph{model\/} (i.e. complete simply connected) surface of constant curvature $\kappa$, $D(\kappa)\/$ be its \emph{diameter}, that is, $D(\kappa)\, =\, \infty\/$ if $\kappa\, \le\, 0\/$ and $D(\kappa)\, =\, \pi/\sqrt{\kappa}\/$ if $\kappa \, >\, 0$. The metric of $M_\kappa\/$ is denoted as $d_\kappa$.

A \emph{triangle\/} in $X\/$ is the union of three mimimizing segments $c_i\,\colon [a_i, b_i] \, \to\, X$, $i\, =\, \overline{1,3}$, called \emph{sides\/} of triangle, pairwise connecting three points $x_i, i\, =\, \overline{1,3}\/$ which we call \emph{vertices\/} of triangle.

For a collection $\Delta := (x_1, x_2, x_3)\/$ of points in $X\/$ the \emph{comparison triangle\/} $\overline{\Delta}\subset M_\kappa\/$ has vertices $\overline{x}_1, \overline{x}_2\/$ and $\overline{x}_3\/$ such that $d_\kappa(\overline{x}_i, \overline{x}_j)\, =\, d(x_i, x_j)$, $i,j=\overline{1,3}$. A comparison triangle exists and is unique up to isometry if the perimeter 
\[P(\Delta) := d(x_1, x_2) \, +\, d(x_2, x_3) \, +\, d(x_1, x_3) \, <\, 2 D(\kappa)\, .\]
If $c_i\/$ are sides of $\Delta$, we denote $\overline{c}_i\/$ sides of $\overline{\Delta}$.

A point $\overline{x}\in \overline{\Delta}\/$ corresponds to a point $x\/$ of triangle $\Delta\/$ if there is some  $i\/$ and  some $t_i\in [a_i, b_i]\/$ with  $x=c_i(t_i)\/$ and $\overline{x}=\overline{c}_i(t_i)$.

Let  $\kappa\in \R$.  A complete metric  space $X\/$ is called  a \emph {$CAT(\kappa)$-space\/} if\\
{(i)} Every two points $x, y\in X\/$ with $d(x, y)\, <\, D(\kappa)\/$ are connected  by a     minimizing;\newline
{(ii)} For any triangle  $\Delta\/$ in $X\/$ with perimeter less than $2D(\kappa)\/$ and any two points $x,y\in \Delta$, the inequality  $d(x,y)\, \le\, d(\overline{x},\overline{y})\/$ holds,  where  $\overline{x}\/$ and $\overline{y}\/$ are the points of the comparison triangle $\overline{\Delta}\subset M_\kappa\/$  corresponding to $x\/$ and $y\/$  respectively.

Every geodesic $c\,\colon \R \, \to\, X\/$ in $CAT(0)$-space $(X,d)\/$ has two opposite directions defined by its rays $c|_{(-\infty, 0]}\/$ and $c|_{[0, +\infty)}$. Geodesic rays $c_i\,\colon [0, \infty) \, \to\, X, i=\overline{1,2}\/$  are called \emph{asymptotic\/} if Hausdorff distance between them is finite: $d_H(c_1, c_2) \, <\, \infty$. A relation \emph{to be asymptotic\/} on a set of all rays in $X\/$ is equivalence. Equivalence classes called \emph{ideal points\/}  (or \emph{points at infinity}) form \emph{geometric boundary\/}  (\emph{boundary at infinity}) $\partial_\infty X\/$ of $X$. Complete geodesics $c_1\/$ and $c_2\/$ are \emph{asymptotic\/} if they contain asymptotic rays.

Complete geodesics $c_1\/$ and $c_2\/$ are called \emph{parallel\/} if their Hausdorff distance ${\rm Hd}(c_1, c_2)\/$ is finite. In this case they are asymptotic in both directions and bound a \emph{flat strip}, i.e. subset isometric to a strip in Euclidean plane.

\begin{definition}
We say that a geodesic $c\/$ \emph{virtually bounds a flat strip\/} if there is a finite sequence $c\, =\, c_0, c_1, \dots, c_n\/$ of geodesics such that for every $1\, \le\, i\, \le\, n\/$ geodesics $c_{i-1}\/$ and $c_i\/$ are asymptotic and $c_n\/$ is a boundary geodesic of some flat strip. In particular every line bounding the flat strip also bounds it virtually. Otherwize we say that $c\/$ \emph{has strictly rank one}.
\end{definition}

\subsection{Spaces of directions}

A \emph{pseudo-metric\/} on a set $\Sigma\/$ is a function $d\,\colon  \Sigma\times \Sigma\, \rightarrow\, [0,   \infty)\/$ that is symmetric and satisfies the triangle inequality. If $(\Sigma,d)\/$ is a  pseudo-metric space, then we get a metric space $(\Sigma^\ast,  d^\ast)\/$ by letting $\Sigma^\ast\/$ be the set of maximal zero
diameter subsets and  setting  $d^\ast(S_1,S_2):=d(s_1,s_2)\/$ for any $s_i\in S_i$. $(X^\ast, d^\ast)\/$ is called the metric space \emph{associated\/} to the pseudo-metric $d$.
 
Given two geodesic segments $c=[x\eta]\/$ and $d=[x\zeta]\/$ with common vertex $x\/$ in $CAT(\kappa)$-space $X$, $\kappa \in \R\/$ with lengthes less than $D(\kappa)$, we have a well-defined function $\angle_x(c,d)\, =\, \lim\limits_{t \, \to\, 0} \widetilde {\angle_x}(c(t),d(t))$, where $\widetilde {\angle_x}(y,z)\/$ denotes \emph{the comparison angle}, i.e. angle at the vertex $\overline x\/$ of comparison triangle $\overline x\,\overline y\,\overline z\/$ at the Euclidean plane for triangle $xyz$. The function $\angle_x\/$ defines a pseudo-metric on a set of all geodesic segments begining from $x$. The metric space associated  to the pseudo-metric  $\angle_x\/$ is denoted by $\Sigma^\ast_x X$. Its metric completion $\Sigma_x X\/$ is called \emph{the space of directions\/} at $x$. If a space $X\/$ is geodesically complete, then we have $\Sigma^\ast_x X\, =\, \Sigma_x X\/$ for any $x \in X$.

\begin{Th}[Nikolaev \cite{N}\label{t8}]
Let $X\/$ be a $CAT(\kappa)$-space and $p\in X$.   Then $\Sigma_p X\/$ is a $CAT(1)$-space.
\end{Th}

\subsection{Diagonal tube and satellite relations} \label{satel}

From now on we assume that $(X,d)\/$ is a locally compact geodesically complete connected at infinity $CAT(0)$-space (\emph{Hadamard space}) and $\overline{X}\, =\, X \cup \partial_\infty X\/$ be its geometric closure. The metric $d\/$ of $CAT(0)$-space is convex function in following sence.
For every two geodesics $c_1,c_2\,\colon \R \, \to\, X\/$  the function $d(c_1(t), c_2(t))\/$ is convex.

For two points $y,z \in \overline{X}\/$ notation $[yz]\/$ means:
\begin{itemize}
\item connecting them geodesic segment, if both lies in $X$, or
\item connecting them geodesic ray if one of them lies in $X$, and another (order is ignored) --- in $\partial_\infty X$, or
\item any complete geodesic connecting  $y\/$ and $z\/$ in the case when $y,z \in \partial_\infty X$, if such geodesic does exist. When the geodesic is not unique, we will detect it by additional features.
\end{itemize}
We will identify every geodesic as a map of real interval into $X\/$ with its image and use the same notation.

For  a number $r\, >\,0\/$\, $V\/$ will be \emph{the diagonal tube\/} corresponding it. Recall that we have assumed $r=1$, hence
\[V := \{(x,y) \in X\times X|\quad |xy|\, \le\, 1\}\, .\]

When we say that diagonal tube $V\/$ defines metric $d$, this means that every metric $d'\/$ on $X\/$ satisfying conditions of theorem \ref{tube} with the same diagonal tube $V\/$ coincides with $d$. We will consider $d\/$ as initial metric and $d'\/$ as trial metric for which we need to show that $d'\, =\, d$. Below $d'\/$ will always be trial metric satisfying conditions of theorem \ref{tube} and with the same diagonal tube $V$.

Also, sometimes we will use the terminology, such as following: "diagonal tube defines ..." or "$V\/$ allows to recover..." etc. Opinions of such type will mean that some object, value or some property of object is the same for metric $d\/$ and any trial metric $d'$. In particular this holds when mensioned object admits description in terms of $V$.

First, we note that  relations  $nV$, ${\rm int}(nV)\/$ and $\partial(nV)\/$ are defined for all $n \in \N\/$ simultaneously with $V\/$ by equalities
\[nV := \{(x,y)\in X\times X |\quad |xy|\, \le\, n\}\]
\[\partial(nV) := \{(x,y) \in X\times X |\quad |xy|\, =\, n \}\] and 
\[{\rm Int} (nV) := nV \setminus \partial (nV)\, .\]

\begin{Le}\label{poli}
If two metrics $d\/$ and $d'\/$ on $X\/$ as in theorem \ref{tube} have common diagonal tube $V$, then they have also common relations $nV$, $\partial(nV)\/$ and ${\rm Int}(nV)\/$ for all $n \in \N$.
\end{Le}

{\sloppy
\begin{proof} Pair $(x,y) \in 2V\/$ iff there exists a point $z \in X\/$ such that $(x,z), (z,y) \in V$. Inductively, $(x, y) \in nV\/$ iff there exists $z \in X\/$ such that $(x,z)\in(n-\nolinebreak 1)V\/$ and $(z,y) \in V$. We have $(x, y) \in \partial (2V)\/$ iff there exists unique point $z\/$ with $(x,z), (z,y) \in V$. In this case $(x,z)\/$ and $(z,y)\in \partial V\/$ and $z\/$ is a midpoint of segment $[xy]$. Moreover, pair $(x,z) \in \partial V\/$ iff there exists a pair $(x,y) \in \partial (2V)\/$ such that $z\/$ is a midpoint of $[xy]$. Inductively one may define relations $\partial (nV)$. Relations $\inter (nV)\/$ are by definition \[\inter (nV)\, =\, nV \setminus \partial (nV)\, .\]
\end{proof}

}
Consequently metrics $d\/$ and $d'\/$ with common $V\/$ have common open, closed balls and spheres with integer radii.

\begin{Le}\label{inter}
If two metrics $d\/$ and $d'\/$ on $X\/$ as in theorem \ref{tube} have common any one of three relations $V$, $\partial V\/$ and $\inter V$, then they have common another two relations.
\end{Le}

\begin{proof} It was shown in Lemma \ref{poli} that $V\/$ defines $\partial V\/$ and ${\rm Int}(V)$.  

Assume that $d\/$ and $d'\/$ have common boundary relation $\partial V$. Arguments similar to those used in the proof of Lemma \ref{poli} show that they have common diagonal tube $2V\/$ as well. In particular, metrics $d\/$ and $d'\/$ have common all 
balls $\overline{B}(x, 2n)\/$ of even radii and consequently all spheres $S(x, 1)\/$ and $S(x, 2n)\/$ with $x \in X$.
Then by the connectedness at infinity of $X\/$ we have $(x, y) \in {\rm int} V\/$ iff $S(x,1) \cap S(y, 1) \ne \emptyset\/$ and $S(x,2) \cap S(y, 1)\, =\, \emptyset$. This defines ${\rm Int} V\/$ by already defined relations. $V\/$ is defined as $V\, =\, {\rm Int} V \cup \partial V$. 

Let now $d\/$ and $d'\/$ have common ${\rm Int} V$. Then for all $n \in \N\/$ relations ${\rm Int}(nV)\/$ are also common for $d\/$ and $d'$. Hence from geodesic completeness of $X\/$ we have that $(x,y) \in \partial V\/$ iff $(x,y) \notin {\rm Int} V\/$ and $B(y,1) \subset B(x,2)$. Indeed, if $1 \, <\, d(x,y) \, <\, 2$, then we may include geodesic segment $[xy]\/$ into segment $[xz]\/$ of length $d(x,z)\, =\, \frac 12 (3 \, +\, d(x,y))$, and $z \in B(y,1) \setminus B(x,2)$. This defines $\partial V\/$ and hence $V\/$ in terms of ${\rm Int} V$. 
\end{proof}

Consequently for the proof of theorem \ref{tube} we need only to prove that the metric $d\/$ may be recovered from diagonal tube $V\/$ itself.

\subsection{$r$-sequences.}

The main tool of \cite{Berestovskii}, called $r$-sequence was defined in terms of relations above. We give another definition here.
\emph{$r$-sequence\/} in  $X\/$ is by definition a homothety $\Z \, \to\, X\/$ with coefficient $r$, that is an integral parametrized sequence $\{ x_z\}_{z \in \Z}\subset X\/$ of points in $X\/$, such that for all $z_1, z_2 \in \Z\/$ equalities $|x_{z_{1}}\, x_{z_{2}}|\, =\, r |z_2\, -\, z_1|\/$ hold.
When we put $r=1$, the last equation becomes written $|x_{z_{1}}\, x_{z_{2}}|\, =\, |z_2\, -\, z_1|\/$ and we consider $r$-sequence as isometric map $\Z \, \to\, X$.

It immediately follows from lemmas \ref{poli} and \ref{inter} and geometry of space $(X, d)\/$ that:
\begin{itemize}
\item for every $r$-sequence there is unique containing it geodesic in $(X,d)\/$ and geodesic in $(X, d')$. A priori these geodesics may be different. If they coincide, we will say that the incidence relation on this geodesic is detected by $V$.
\item for any trial metric $d'$, sequence $\{x_z \}_{z \in \Z}\/$ is $r$-sequence with respect to the metric $d\/$ iff it is $r$-sequence with respect to $d'$, and
\item $V\/$ allows us to reveal, given two $r$-sequences $\{x_z \}_{z \in \Z}\/$ and $\{y_z \}_{z \in \Z}\/$ whether geodesics in $(X,d)\/$ containing them are asymptotic in any direction.
\end{itemize}

The segment of $r$-sequence $\{ x_z\}_{z \in \Z}\subset X\/$ between points $x_z\/$ and $x_{z\, +\,k}\/$ will be denoted as $[x_z, \dots , x_{z\, +\,k}]_r$, ideal points, defined by $r$-sequence $\{x_z \}_{z \in \Z}$, --- as $x_{+\infty}\/$ and $x_{-\infty}$.

\subsection{Metric transfer}\label{trans}
We will use an effective trick proposed in \cite{Berestovskii}, --- a map $R_{c_{1}c_{2}}\/$ from a geodesic to its asymptotic one by means of \emph{Busemann functions}. Let $c_i \,\colon \R \, \to\, X, i\, =\, \overline{1, 2}\/$ be two asymptotic geodesics such that $c_1(+\infty)\, =\, c_2(+\infty)$. Let $\beta_\xi\/$ be some Busemann function corresponding to an ideal point $\xi\, =\, c_i(+\infty)$:
\begin{equation}\label{hs}
\beta_\xi (x) := \lim\limits_{t \, \to\, \infty}(|x\, c(t)|\, -\, t).
\end{equation}
Level sets of function $\beta_\xi\/$ are called \emph{horospheres}, sublevels --- \emph{horoballs}. We will distinguish open horoballs defined by strict inequality $\beta_\xi(x) \, <\, b_0\/$ and horoballs $\beta_\xi(x)\, \le\, b_0$.

The horosphere centered in the end $\xi := c(+\infty)\in \partial_\infty X\/$ of geodesic $c\,\colon\R \, \to\, X\/$ containing a point $z \in X$, i.e. level set \[\{x \in X| \, \beta_\xi(x)\, =\, \beta_\xi(z)\}\] of Busemann function $\beta_\xi(x)\/$ \eqref{hs} will be denoted as $\mathcal{HS}_{\xi, z}$, corresponding closed horoball as $\mathcal{HB}_{\xi,z}$. An open horoball is $hb_{\xi,z} := \mathcal{HB}_{\xi,z}\setminus \mathcal{HS}_{\xi,z}$.

Horoballs and horospheres admits following definitions in terms of ${\rm Int} V\/$ (and hence in terms of $V$). Let $\{x_z\}_{z \in \Z}\/$ be $r$-sequence with $\xi := x_{+\infty}\, =\, c_i(+\infty)\/$ and given point $x_0$. Then open horoball $hb_{\xi, x_0}\/$ is
\[hb_{\xi, x_0}\, =\, \bigcup\limits_{n=1}^{\infty} B(x_n, n)\, .\]
For a point $y \in X\/$ we have $y \in \mathcal{HS}_{\xi,x_0}\/$ iff $y \notin hb_{\xi, x_0}\/$ and $B(y,1) \subset hb_{\xi,x_{-1}}$. Finally, $\mathcal{HB}_{\xi,x_0}\, =\, hb_{\xi,x_0} \cup \mathcal{HS}_{\xi,x_0}$.

There exist length parametrisations of geodesics $c_i$, such that 
\begin{equation}\label{parametr}
\forall t \in \R \quad \beta_\xi (c_1(t))\, =\, \beta_\xi (c_2(t)).
\end{equation} 

We set $R_{c_{1}c_{2}} (c_1(t)) := c_2(t)$. A map $R_{c_{1}c_{2}}\/$ defined above is independent on a choise of Busemann function $\beta_\xi\/$ centered in $c_i(+\infty)\/$ and on parametrisations with the property \eqref{parametr}. 

It was shown in  \cite{Berestovskii} that whenever the relation of incidence of points in lines  $c_i\/$ is detected, then relation $V\/$ defines every map  $R_{c_{1}c_{2}}\/$ in metric space $(X, d)\/$ itself. This means that whenever geodesics $c_1\/$ and $c_2\/$ are also geodesics with respect to trial metric $d'$, then the image of any point $c_1(t)\/$ under the map $R_{c_{1}c_{2}}\/$ is independent on the choise of metric $d\/$ or $d'\/$.

\begin{definition}
We say that two geodesics $a\/$ and $b\/$ are \emph{connected by the asymptotic chain\/} if there is a finite sequence 
\begin{equation}\label{a1n}
a\, =\, a_0, a_1, \dots, a_n=b
\end{equation}
of geodesics such that for every $1\, \le\, i\, \le\, n\/$ geodesics $a_{i-1}\/$ and $a_i\/$ are asymptotic in some direction.
\end{definition}
\begin{example}
If geodesic $a\/$ virtually bounds flat strip, it is connected by asymptotic chain with some geodesic $b\/$ lying in the boundary of flat strip.
\end{example}

\begin{Le}\label{asymchain}
Let geodesics $a\/$ and $b\/$ are connected by the asymptotic chain {\rm \eqref{a1n}} such that for all $i\, =\, \overline{0,n}\/$ relation of incidence in geodesics $a_i\/$ is detected by $V$. Then if distances $d\/$ and $d'\/$ coincide along $a$, i.e. for every pair $t_1, t_2\in \R\/$ equality $d(a(t_1), a(t_2))\, =\, d'(a(t_1), a(t_2)\/$ holds, then metrics coincide along $b\/$ as well.
\end{Le}

\begin{proof}
The claim may be proved by a multiple transfer of a metric from a geodesic to its asymptotic one in the asymptotic chain connecting $a\/$ and $b$.
\end{proof}

\subsection{Parallel equivalence}\label{paral}

\begin{definition}
Let $(\{x_z\}_{z \in \Z},\, \{y_z\}_{z \in \Z})\/$ be a pair of $r$-sequences in $X$. We say that they are  \emph{parallel-equivalent}, iff Hausdorff distance $d_H\/$ between sets $\{x_z\}_{z \in \Z}\/$ and $\{y_z\}_{z \in \Z}\/$ is finite: 
\[d_H\left( \{x_z\}_{z \in \Z},\, \{y_z\}_{z \in \Z}\right) \, <\, +\infty\, .\]
\end{definition}
We will use the notation $\{x_z\}_{z \in \Z}\, \parallel \{y_z\}_{z \in \Z}\/$ for this relation. It is obviously a really equivalence on a set of all $r$-sequences in $X$. In fact every $r$-sequence belongs to a unique geodesic in $X$, and  $\{x_z\}_{z \in \Z} \parallel \{y_z\}_{z \in \Z}\/$ iff  geodesics containing $\{x_z\}_{z \in \Z}\/$ and $\{y_z\}_{z \in \Z}\/$  are parallel or coinside.

\begin{Le}
Property "to be parallel-equivalent" can be revealed by $V$.
\end{Le}

{\sloppy
\begin{proof} $\{x_z\}_{z \in \Z}\parallel\{y_z\}_{z \in \Z}\/$ iff there  exists $k\in \N\/$ such that for all $z\in\Z\/$ $(x_z, \, y_z) \in kV$. This completes the proof.
\end{proof} 

}

\begin{definition}
$r$-sequence is called $r$-sequence \emph{of rank\/} 1, if there exists unique geodesic $c\/$ in $X\/$ containing all $r$-sequences parallel-equivalent to it. $c\/$ has not different parallel geodesics in $X\/$ and does not bound a flat strip in this case. It was shown in \cite{Berestovskii} that the incidence relation is detected by the diagonal tube in this case. In opposite case we say that $r$-sequense is \emph{of higher rank\/} or has \emph{rank over\/} 1. Geodesic containing $r$-sequence of higher rank bounds some flat strip. 
\end{definition}

\begin{Prop} \label{rank}
A property of an arbitrary $r$-sequence $(\{x_z\}_{z \in \Z}\/$ to be of rank 1 or of higher rank is determined by relation $V$.
\end{Prop}
\begin{proof}
$\{x_z\}_{z \in \Z}\/$ is of higher rank iff there exists $r$-sequence $\{ y_z \}_{z \in \Z}\/$ parallel-equivalent to it and not equal to $\{ x_{z\pm 1} \}_{z \in \Z}$, such that $| x_0\, y_0 |\, =\, 1$. Since the relation $\partial V\/$ is determined by $V\/$ (cf Lemma \ref{poli}), the proposition is proved.
\end{proof}

The plan of further consideration is following. We will examine every single geodesic in accordance with its rank. It will be shown that in any case the incidence relation and the metric of a geodesic may be recovered from a diagonal tube $V$. First we consider a case of geodesic bounding a flat strip.

\section{The geodesic bounding flat strip}
\subsection{Splitting of parallel-equivalence class}\label{split}

Let $c \,\colon (-\infty , +\infty) \, \to\, X\/$ be a geodesic connecting ideal points $\xi_- := c(-\infty)\/$ and $\xi_+ := c(+\infty)\/$ and passing through the point $x_0\, =\, c(0)$. For $x \in X\/$ and $\xi, \eta \in \partial_\infty X$, angle at point $x\/$ between rays $[x \, \xi]\/$ and  $[x \, \eta]\/$ is denoted as $\angle_x (\xi, \, \eta)$.

The following lemma is well-known (cf.  \cite[Lemma 5]{Okun} for example).

\begin{Le}\label{secheniya}
Set \[C := \{ y \in X |\, \angle_y (\xi_-,\, \xi_+)\, =\, \pi\}\, .\] Then
\begin{enumerate}
\item $\mathcal{HS}_{\xi_-, \, x_0} \cap \mathcal{HS}_{\xi_+, \, x_0}\, =\,  \mathcal{HB}_{\xi_-, \, x_0} \cap \mathcal{HB}_{\xi_+, \, x_0}$,
\item $C\/$ is a union of all geodesics parallel to $c$,
\item $C\/$ is a closed convex set in  $X\/$ and
\item $C\/$ splits as a product $C\, =\, C\,' \times c$, where $C\,' := \mathcal{HS}_{\xi_-, \, x_0} \cap \mathcal{HS}_{\xi_+, \, x_0}\/$ is a closed convex subset.
\end{enumerate}
\end{Le}

Hence the set $C\/$ of all points of all $r$-sequences parallel-equivalent to given  $r$-sequence  $\{x_z\}_{z \in \Z} \subset X\/$ splits as $C\, =\, C\,' \times c$. The set $C\,'\/$ has more than one point iff the $r$-sequence $\{x_z\}_{z \in \Z} \subset X\/$ is of higher rank.

\begin{remark}
Given arbitrary $r$-sequence $\{x_z'\}_{z \in \Z}\/$ parallel-equivalent to $r$-sequence $\{x_z\}_{z \in \Z} \subset X$, denote as $C\,'(x_z')\/$ the fiber of a product $C=C\,'\times~c\/$ containing a point $x_z'$. Since horospheres  $\mathcal{HS}_{\xi_-, \, x_0}\/$ and $\mathcal{HS}_{\xi_+, \, x_0}\/$ may be defined in terms of $V$, the \emph{horisontal structure\/} of splitting is also determined by $V$. In other words any trial metric $d'\/$ gives rise to the same set of fibers of type $C\,'(x)$, which we will call \emph{horisontal sections\/} of $C$.
\end{remark}

As well, diagonal tube $V\/$ allows to restore the order of horisontal sections of  $C\/$ as following. We say that a section $C\,'(x)\/$ of the set $C\/$ \emph{lays below\/} of a section $C\,'(y)\/$ if a ray $[x \xi_+]\/$ intersects  $C\,'(y)$. We have $\mathcal{HB}_{\xi_-, \, x} \cap \mathcal{HB}_{\xi_+, \, y}\, =\, \emptyset\/$ in this case. It is left to restore distances between horisontal sections of $C\/$ and incidence of points of $c\/$ for completing the case.
 
\subsection{Tapes}

\begin{definition}
We say that $4p \quad (p \in \N)\/$ pairwise parallel-equivalent  $r$-sequences 
\begin{equation}\label{baza}
\{x_{i,\, j;\, z}\}_{z \in \Z},\, i= \overline{0,\, 3}, j= \overline{1,\, p}
\end{equation}
form \emph{$p$-tape}, if following $4p\, +\,4\/$ points 
\[
\begin{array}{l}
x_{i,\, 1,\, 0}, \dots, x_{i,\, p,\, 0}, \, i\, =\, \overline{0,\, 3} \\
x_{0,\, 1,\, 2p-1},\, x_{2,\, p,\, 1-2p},\, x_{3,\, p-1,\, 1-2p}, x_{3,\, p,\, 1-2p}
\end{array}
\]
in addition satisfy the system of $2p\/$ relations
\begin{equation}
\label{tesyomka}
\begin{cases}
		[x_{0,\, 1,\, 0},\, x_{1,\, 1,\, 0},\, x_{2,\, 1,\, 0},\, x_{3,\, 1,\, 0}]_r \\
		\qquad\qquad\qquad\dots \\
		\left[x_{0,\, p,\, 0},\, x_{1,\, p,\, 0},\, x_{2,\, p,\, 0},\, x_{3,\, p,\, 0}\right]_r \\
		\left[x_{0,\, 2,\, 0},\, x_{1,\, 1,\, 0},\, x_{2,\, p,\, 1-2p},\, x_{3,\, p-2,\, 1-2p}\right]_r\\
		\left[x_{0,\, 3,\, 0},\, x_{1,\, 2,\, 0},\, x_{2,\, 1,\, 0},\, x_{3,\, p-1,\, 1-2p}\right]_r\\
		\left[x_{0,\, 4,\, 0},\, x_{1,\, 3,\, 0},\, x_{2,\, 2,\, 0},\, x_{3,\, 1,\, 0}\right]_r\\
		\qquad\qquad\qquad\dots\\
		\left[x_{0,\, p,\, 0},\, x_{1,\, p-1,\, 0},\, x_{2,\, p-2,\, 0},\, x_{3,\, p-3,\, 0}\right]_r\\
		\left[x_{0,\, 1,\, 2p-1},\, x_{1,\, p,\, 0},\, x_{2,\, p-1,\, 0},\, x_{3,\, p-2,\, 0}\right]_r
\end{cases}
\end{equation}
presenting segments of $r$-sequences.
\end{definition}

Since $r$-sequences admit definition in terms of $V$, notion of $p$-tape is independent on the choise of metric $d\/$ or trial metric $d'$. A fragment of a $p$-tape is shown at figure  \ref{tesma}. Here only a part of points of $r$-sequences forming $p$-tape is picked out. The main idea is that a segment between points  $x_{0,\, 1,\, 0}\/$ and  $x_{0,\, 1,\, 2p-1}\/$ of the same $r$-sequence contains $2(p-1)\/$ points of this $r$-sequence dividing corresponding geodesical segment by $2p-1\/$ equal parts, and simultaneously marked segment contains $p-1\/$ points of kind $x_{0,\, j, \, 0} \quad j\, =\, \overline{2,\, p}$, dividing it  by $p\/$ equal segments. Hence all points of  $r$-sequences $\{ x_{0, \, j,\, z} \}_{z \in \Z}\/$ at segment $[x_{0,\, 1,\, 0}\, x_{0,\, 1,\, 2p-1}]\/$ divide it by $p(2p-1)\/$ equal parts, and in particular segment $[x_{0, \, 1,\, 0}\, x_{0,\, 1,\, 1}]\/$ is divided by points $x_{0,\, p,\, 3-2p},\, x_{0,\, p-1,\, 5-2p}\,,\, \dots ,\, x_{0,\, 2,\, -1}\/$ by $p\/$ equal segments.

\begin{figure}
	\begin{picture}(320,60)
		\multiput(0,10)(80,0){3}{\line(5,1){120}}
		\multiput(40,34)(80,0){3}{\line(5,-1){120}}
		\put(0,26){\line(5,-1){80}}
		\put(0,26){\line(5,1){40}}
		\put(240,10){\line(5,1){80}}
		\put(280,34){\line(5,-1){40}}
		\multiput(0,10)(80,0){5}{\circle*{3}}
		\multiput(40,18)(80,0){4}{\circle*{3}}
		\multiput(0,26)(80,0){5}{\circle*{3}}
		\multiput(40,34)(80,0){4}{\circle*{3}}
		\put(0,0){$x_{0,\, 1,\, 0}$}
		\put(110,40){$x_{3,\, 1,\, 0}$}
		\put(225,0){$x_{0,\, p,\, 0}$}
		\put(295,0){$x_{0,\, 1,\, 2p-1}$}
		\put(18,40){$x_{3,\, p,\, 1-2p}$}
		\put(70,0){$x_{0,\, 2,\, 0}$}
		\put(70,35){$x_{2,\, 1,\, 0}$}
		\put(30,5){$x_{1,\, 1,\, 0}$}
		\put(260,40){$x_{3,\, p-1,\, 0}$}
		\put(265,5){$x_{1,\, p,\, 0}$}
	\end{picture}
\caption{$p$-tape}
\label{tesma}
\end{figure}

Now the plan of restoring the metric of geodesic of higher rank is following.
Given $r$-sequence  $\{x_z \}_{z \in \Z}\/$ of higher rank, geodesic $c\,\colon(-\infty,+\infty)\, \to\,~X\/$ containing it spans a flat strip (not unique in general). For sufficiently great $P\/$ for any $p \, >\, P$, \quad $p$-tape defined by a family of $r$-sequences of type \eqref{baza} with $x_{0,\, 1,\, z}\, =\, x_z\/$ for all $z\in \Z$, does exist. Points of all such $r$-sequences  $x_{0,\, j,\, z}\/$ for every possible  $p$-tapes  with $p \, >\, P\/$ cover a set of rational points on considered geodesic.

As it was mensioned earlier, horospheres and horoballs of $CAT(0)$-space $X\/$ define horisontal sections of a set $C$. Hence the metric of a term $\R\/$ in the splitting $C\, =\, C' \times \R\/$ of the set $C\/$ above may be recovered from the relation $V$.

It remains to show, that the incidence of points on $c\/$ is restorable by $V$. It is already done for rational points. The limiting procedure using balls is possible for irrationals.

\subsection{Recovery of the metric of geodesic bounding flat strip}

\begin{Le}
Every $p$-tape in $X\/$ is contained in a flat strip.
\end{Le}

\begin{proof} Let $p$-tape be defined by  $r$-sequences \eqref{baza}. Consider a flat strip $F\/$ spanned by parallel geodesics $c\/$ and $c'$, containing correspondingly $r$-sequences $\{ x_{0,\, 1,\, z}\}_{z \in \Z}\/$ and $\{ x_{3,\, 1,\, z}\}_{z \in \Z}$. $r$-sequences $\{ x_{1,\,1,\, z }\}_{z \in \Z}\/$ and $\{ x_{2,\,1,\, z }\}_{z \in \Z}\/$ are also contained in $F$. {\sloppy 

}

We will show that $F\/$ really contains all $r$-sequences in  \eqref{baza}. Consider points 
\begin{equation}\label{6-points}
x_{0,\, 2,\, 0}, \quad x_{1,\, 1,\, 0}, \quad x_{2,\, 1,\, 0}, \quad x_{3,\, 1,\, 0}, \quad x_{2,\, 2,\, 0}, \quad x_{1,\, 2,\, 0}. 
\end{equation}

As a corollary of \eqref{tesyomka} we have relations:
\[[x_{0,\, 2,\, 0}, \, x_{1,\, 2,\, 0}, \, x_{2,\, 2,\, 0}]_r\]
and
\[[x_{1,\, 1,\, 0}, \, x_{2,\, 1,\, 0},\, x_{3,\, 1,\, 0}]_r\, .\]

Because of convexity of metric in a space of nonpositive curvature the function
\[d_{2,\,  1}(t)\, =\, |\, x_{t,\, 2,\, 0}\, x_{t + 1,\, 1,\, 0}\, |\, ,\]  
where $x_{t, \, j,\, 0}\, =\, c_j(t)$,   $c_j\,\colon [0, 2r] \, \to\, X, \quad j\, =\, \overline{1,2}\/$ --- geodesic segments connecting $x_{0,\, j,\, 0}\/$ with $x_{1, \, j,\, 0}$, is a convex function. But \[d_{2,\, 1}(0)\, =\, d_{2,\, 1}(1)\, =\, d_{2,\, 1}(2)\, =\, 1\, ,\]
hence $d_{2,\, 1}(t)\, =\, 1\/$ for all $t \in [0,\, 2]$, and points \eqref{6-points} belong to a flat parallelogram
\begin{equation}\label{parallelogramm-1}
x_{0,\, 2,\, 0} \, x_{1,\, 1,\, 0} \, x_{3,\, 1,\, 0}, \, x_{2,\, 2,\, 0}
\end{equation} 
isometrically embedded into $X$.

Similarly six points
\[x_{2,\, 1,\, 0}, \quad x_{3,\, 1,\, 0}, \quad x_{2,\, 2,\, 0}, \quad x_{1,\, 3,\, 0}, \quad x_{0,\, 3,\, 0}, \quad x_{1,\, 2,\, 0}\]
belong to isometrically embedded into $X\/$ flat parallelogram
\begin{equation}\label{parallelogramm-2}
x_{2,\, 1,\, 0} \, x_{3,\, 1,\, 0} \, x_{1,\, 3,\, 0} \, x_{0,\, 3,\, 0}.
\end{equation}

Parallelograms \eqref{parallelogramm-1} and \eqref{parallelogramm-2} have the common part  --- parallelogram
\[x_{2,\, 1,\, 0} \, x_{3,\, 1,\, 0} \, x_{2,\, 2,\, 0} \, x_{1,\, 2,\, 0}\, ,\]
and hence their union is nonconvex flat hexagon
\[x_{0, \, 2, \, 0} \, x_{1,\, 1,\, 0} \, x_{3,\, 1,\, 0} \, x_{1,\, 3,\, 0} \, x_{0,\, 3,\, 0} \, x_{1, \, 2, \, 0}\, .\]

It is isometric to a nonconvex hexagon in Euclidean plane, constructed by two equal parallelograms. 
Continuing so on we get a flat poligon
\begin{equation}\label{poligon}
\begin{array}{c}
P\, =\, x_{0, \, 2, \, 0} \, x_{1,\, 1,\, 0} \, x_{3,\, 1,\, 0} \, x_{2,\, 2,\, 0} \, x_{3, \, 2, \, 0} \, \dots \\
\dots x_{3, \, p-1, \, 0} \, x_{2, \, p, \, 0} \, x_{3, \, p, \, 0} \, x_{2, \, 1, \, 2p-1}\, x_{0, 1, 2p-1} \, x_{0, \, p, \, 0} \,
x_{1, \, p-1, \, 0} \, \dots \\
\dots x_{1,\, 3,\, 0} \, x_{0,\, 3,\, 0} \, x_{1, \, 2, \, 0},
\end{array}
\end{equation}
(see Figure \ref{tesma}) which is isometric to a $(4p-2)$-gon constructed from $2p-3\/$ parallelograms, consecutively intersecting each other. 

$P\/$ containes whole geodesic segments  $x_{1,\, 1, \, 0}\, x_{1, \, 1, \, 2p-1}\/$ and $x_{2,\, 1, \, 0}\, x_{2, \, 1, \, 2p-1}$. Hence intersection of $P\/$ with a flat strip $F\/$ contains parallelogram
\[x_{1,\, 1, \, 0}\, x_{2,\, 1, \, 0}\, x_{2, \, 1, \, 2p-1} x_{1, \, 1, \, 2p-1}\, .\] Furthemore $P\/$ is contained in parallelogram $x_{0,\, 1, \, 0}\, x_{3,\, 1, \, 0}\, x_{3, \, 1, \, 2p-1} x_{0, \, 1, \, 2p-1}\/$ and hence in $F$. Since there are points of all $r$-sequences \eqref{baza} in $P$, we may say that the whole tape is contained in $F$. 
\end{proof}

\emph{The width\/} of $p$-tape  \eqref{tesyomka} in $X\/$ is a distance between parallel geodesics including $r$-sequences $\{x_{0, \, j,\, z}\}_{z \in \Z}\/$ and $\{x_{3, \, j,\, z}\}_{z \in \Z}$.

\begin{Le}
The width $s(p)\/$ of $p$-tape equals to \[s(p)\, =\, \frac{3\sqrt{4p-1}}{2p}\, .\]
\end{Le}

\begin{proof}
Since every flat strip in $CAT(0)$-space is isometric to a strip in Euclidean plane, it suffuces to calculate a width of standard Euclidean $p$-tape.
\end{proof}

\begin{Cor}
For every flat strip containing a geodesic  $c$, there exists a number $P \, >\, 0$, such that for all $p \, >\, P\/$ the strip contains a $p$-tape generated by a family of $r$-sequences of type {\rm \eqref{baza}} with $\{ x_{0,\, j,\, z}\} \subset c$.
\end{Cor}

{\sloppy
\begin{Le}\label{rational}
Let the geodesic $c\,\colon (-\infty,\, +\infty ) \, \to\, X\/$ contains $r$-sequence $\{ x_z \}_{z \in \Z}$, $x_0\, =\, c(0)\/$ of higher rank. Then for every rational  $q := \frac mn\/$ and $p := kn$, multiple of  $n$, every $p$-tape for which $x_{0, \, 1,\, z}=x_z\/$ contains a point $c(q)\, =\, x_{0,\, j,\, z'}\/$ with 
\[j := n \, +\, 1\, -\, km'\]
and
\[z' := q' \, +\, 1\, -\, 2p \, +\, 2k(m' \, +\, 1)\, ,\]
where $q\,' := \left[q-\frac 1n\right]\/$ --- integral  and  $\frac{m'}{n} := \left\{ q\, -\, \frac 1n \right\}\/$ --- fractional parts of the number $q-\frac 1n$.
\end{Le}

}
\begin{proof}
By direct calculation for a standard flat $p$-tape.
\end{proof}

This yelds, that for all $p\/$ multiple to $n\/$ every $p$-tape constructed with base points $x_z\/$ of given $r$-sequence $\{x_z\}_{z \in \Z}$, contains points $c(\frac mn)\/$ as elements with defined by $p\/$ multiindex. We are ready now to prove the part of Theorem \ref{tube} in the case we consider.

\begin{Th}\label{parabolic}
Let metric space $(X,d)\/$ and trial metric $d'\/$ on $X\/$ be as in Theorem \ref{tube}.
Assume that $c \,\colon (-\infty , +\infty) \, \to\, X\/$ is a geodesic in metric $d\/$ bounding flat strip in $X$. Then $c\/$ is geodesic bounding flat strip in metric $d'\/$ and for all $t_1, t_2 \in \R\/$ 
\begin{equation}\label{dd1}
d'(c(t_1), c(t_2))\, =\, d(c(t_1), c(t_2)).
\end{equation}
\end{Th}

\begin{proof} 
We may assume that $t_1\, =\,0$.
Let $\{x_z\}_{z \in \Z}\/$ be $r$-sequence lying in $c\/$ such that $x_z\, =\, c(z)$. It is $r$-sequence of higher rank relatively both metrics $d\/$ and $d'\/$ and every $p$-tape in metric $d\/$ based on $\{x_z\}_{z \in \Z}\/$ is $p$-tape in metric $d'\/$ as well. When $t\/$ is rational, equality \eqref{dd1} is a consequence of Lemma \ref{rational}.

Let the number  $t\/$ be irrational and $t_n, t'_n \, \to\, t\/$ be a pair of sequences of its rational approximations from below and from above correspondingly. Without loss of generality we may assume that $t\, >\,0$. Then 
\[c(t)\, =\, \left(\bigcap\limits_{n\, =\, 1}^\infty B(c((t'_n\, -\, 1)), \, 1)\right) \cap  \left(\bigcap\limits_{n\, =\, 1}^\infty B(c((t_n \, +\, 1)), \, 1)\right)\, ,\]
and for metric $D \in \{d, d'\}\/$  we have
\[t_n   \, <\, D(c(0), c(t)) \, <\, t'_n\, ,\]
implying the claim. \end{proof}

\begin{remark}
In other words, Theorem \ref{parabolic} asserts that the diagonal tube $V\/$ allows to recover the metric of geodesic $c\/$ bounding a flat strip.
\end{remark}
\begin{Cor}\label{chain}
$V\/$ recovers a metric of any geodesic which virtually bounds a flat strip in $X$.
\end{Cor}

\begin{proof} 
The claim is a partial case of lemma \ref{asymchain}, since every geodesic in metric $d\/$ does not bounding flat strip is also geodesic in trial metric $d'$.
\end{proof}

\section{Geodesics of rank one}

\subsection{Scissors}

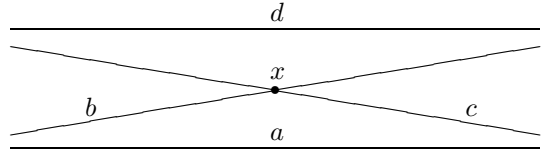
\begin{figure}
\begin{picture}(200,50)
\put(0,0){\line(1,0){200}}
\put(0,5){\line(6,1){200}}
\put(200,5){\line(-6,1){200}}
\put(0,45){\line(1,0){200}}
\put(100,22){\circle*{3}}
\put(98,26){$x$}
\put(98,3){$a$}
\put(98,48){$d$}
\put(172,12){$c$}
\put(28,12){$b$}
\end{picture}
\caption{Scissors $\langle a,b,c,d ;x \rangle$}
\label{scissors}
\end{figure}

\begin{definition}
We say that four complete geodesics $a, b, c, d\,\colon (- \infty , +\infty) \, \to\, X \/$ in $CAT(0)$-space $X\/$ form \emph{scissors centered\/} in $x\in X\/$ if:
\begin{itemize}
\item $a(-\infty)=b(-\infty)$; 
\item $a(+\infty)=c(+\infty)$; 
\item $c(-\infty)=d(-\infty)$;  
\item $b(+\infty)=d(+\infty)$; and
\item $b \cap c\, =\, x$.
\end{itemize}
Such configuration will be denoted as $\langle a,b,c,d ;x \rangle\/$ (figure \ref{scissors}). Geodesics $a\/$ and $d\/$ above are called \emph{bases\/} of scissors. One of them, for example $a$, will be labeled as \emph{lowest base}. Writing $\langle a,b,c,d ;x \rangle\/$ we put the lowest base first. One or both of bases may in general pass throw $x_0$, but when $X\/$ has a property of nonbranching of geodesics, such a situation is excluded.
\end{definition}

\begin{figure}
\begin{picture}(200,70)
\put(0,0){\line(5,1){50}}
\put(50,10){\line(1,0){100}}
\put(150,10){\line(5,-1){50}}
\put(50,10){\line(2,1){100}}
\put(50,60){\line(2,-1){100}}
\put(0,70){\line(5,-1){50}}
\put(50,60){\line(1,0){100}}
\put(150,60){\line(5,1){50}}
\put(100,35){\circle*{3}}
\put(98,38){$x$}
\put(98,13){$a$}
\put(98,63){$d$}
\put(130,22){$c$}
\put(65,22){$b$}
\end{picture}

\vspace{1cm}
\begin{picture}(200,50)
\put(0,0){\line(3,1){90}}
\put(90,30){\line(1,0){20}}
\put(110,30){\line(3,-1){90}}
\put(0,60){\line(3,-1){90}}
\put(110,30){\line(3,1){90}}
\put(100,30){\circle*{3}}
\put(98,33){$x$}
\put(25,3){$a$}
\put(170,3){$a$}
\put(25,55){$d$}
\put(170,55){$d$}
\put(135,12){$c$}
\put(65,42){$c$}
\put(135,42){$b$}
\put(65,12){$b$}
\end{picture}

\caption{Examples of scissors with branching of lines}
\label{scisdeg}
\end{figure}
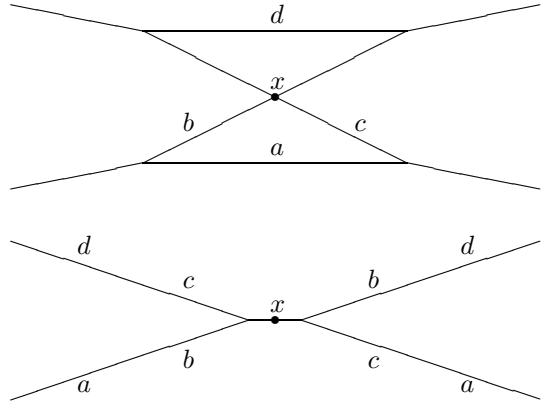

Note that when scissors are given, they are not uniquelly determined by a choice of their lowest base $a\/$ and a center $x\/$  in general: 
we must take into account the possibility of branching of geodesics $b\/$ and $c$. Highest bases $d_1\/$ and $d_2\/$ of scissors $\langle a, b, c, d_1, x \rangle\/$ and   $\langle a, b, c, d_2, x \rangle\/$ with common lines $a$, $b$, $c\/$ may be parallel. Lines forming scissors may be partially or even totally attached. Some examples of scissors with branching of their lines are presented at figure~\ref{scisdeg}.

The main advantage of using scissors consists in a  scissors translation $T$, along their lowest base. It may be described as following. Assume that $R_{ac}\/$ is metric transfer from  geodesic line $a\/$ to the line $c\/$ generated by Busemann function $\beta_{a(+\infty)}\/$ as in subsection \ref{trans}: any point $m \in a\/$ moves to the unique point $m^\prime\, =\, R_{ac}(m) \in c\/$ with $\beta_{a(+\infty)}(m^\prime)\, =\, \beta_{a(+\infty)}(m)$. Similarly maps  $R_{cd}$, $R_{db}\/$ and $R_{ba}\/$ are defined as metric transfers of corresponding lines. All maps above do not depend on the chois of Busemann functions $\beta\/$ in classes defined by corresponding ideal points.
\begin{definition}
\emph{Scissors translation\/} $T\/$ is a composition \[T:=R_{ba}\circ R_{db} \circ R_{cd} \circ R_{ac} \,\colon a \, \to\, a\, .\] Since all maps forming the translation  $T\/$ above are isometric maps, $T\/$ is an isometry of a geodesic $a\/$ preserving its direction. \emph{The displacement\/} of the tanslation $T\/$ i.e. difference $\beta_{a(-\infty)}(T(m))\, -\, \beta_{a(-\infty)}(m)\/$ does not depend on the choice of the point $m\in a$. We denote this value as $\delta T$.
\end{definition}

The term $\delta T\/$ admits following approach. Let $\beta_{a-}$, $\beta_{a+}$, $\beta_{d-}\/$ and $\beta_{d+}\/$ be four Busemann functions such that there exists points $p \in a\/$ and $q \in d$, satisfying equalities $\beta_{a-}(p)\, =\, \beta_{a+}(p)\, =\, 0\/$ and $\beta_{d-}(q)\, =\, \beta_{d+}(q)\, =\, 0$. 

\begin{Th}
\begin{equation}\label{15.02.03-1}
\delta T\, =\, \beta_{a-}(x) \, +\, \beta_{a+}(x) \, +\, \beta_{d-}(x) \, +\, \beta_{d+}(x)\, \ge\, 0,
\end{equation}
where $x\/$ is a center of scissors $\langle a, b, c, d; x \rangle\/$ with translation $T$. Moreover, if none of geodesics $a\/$ or $d\/$ bounds flat strip and if $a\cap d\, =\, \emptyset$, then $\delta T \, >\, 0$
\end{Th}

\begin{proof} First, note that sums  $\beta_{a-}(x) \, +\, \beta_{a+}(x)\/$ and $\beta_{d-}(x) \, +\, \beta_{d+}(x)\/$ are independent on the choice of points $p \in a\/$ and $q \in d$, since,  for example, substituting point $p\/$ by $p' \in a$, one adds constants $\beta_{a-}(p')\/$ and $\beta_{a+}(p')$, where $\beta_{a-}(p')\, =\,\, -\, \beta_{a+}(p')$, to functions $\beta_{a-}\/$ and $\beta_{a+}\/$ correspondingly. We have for $\delta T$:
\[\delta T\, =\, t'\, -\, t\, ,\]
where 
\[a(t')\, =\, T(a(t))\, .\]

If $d(s)\, =\, R_{cd}\circ R_{ac} (a(0))$, then for all $t \in \R\/$ \[d(s\, +\,t)\, =\, R_{cd} \circ R_{ac}(a(t))\, .\] Analogously, if $a(t)\, =\, R_{ba} \circ R_{db}(d(0))$, then for all $s \in \R\/$  \[a(t \, +\, s)\, =\, R_{ba} \circ R_{db}(d(s))\, .\] 

Set $p := a(0)\, =\, R_{ac}^{-1}(x)\/$ and $q := d(0)\, =\, R_{db}^{-1}(x)$. Then \[\beta_{a+}(x)\, =\, \beta_{a+}(p)\, =\, 0\, ,\] \[\beta_{d+}(x)\, =\, \beta_{d+}(q)\, =\, 0\] and 
\[T(p)\, =\, R_{ba} \circ R_{db} \circ R_{cd} (x)\, =\, R_{ba} \circ R_{db} (d(s))\, ,\]
where $s= \beta_{d-}(x)$. We have:
\[T(p)\, =\, a(\beta_{d-}(x) \, +\, t)\, ,\]
with $t := \beta_{a-}(x)$. Hence displacement of a point $p \in a$, and consequently of any point of $a\/$ equals to

\[\delta T\, =\, \beta_{a-}(x) \, +\, \beta_{d-}(x)\, -\, 0\, =\, \beta_{a-}(x) \, +\, \beta_{a+}(x) \, +\, \beta_{d-}(x) \, +\, \beta_{d+}(x)\, .\]

We have $\delta T \, \ge\, 0$, since $x\/$ lies in intersections $hb(a(+\infty), x) \cap hb(a(-\infty), x)\/$ and $hb(d(+\infty), x) \cap hb(d(-\infty), x)$ of horoballs. When additional conditions are cassumed, we may suppose without lost of generality, that $x\notin a$. In this case $\beta_{d-}(x) \, +\, \beta_{d+}(x) \, \ge\, 0\/$ and $\beta_{a-}(x) \, +\, \beta_{a+}(x) \, >\, 0$.

\end{proof}

\subsection{Shadows}\label{Sha}

\begin{definition}
Fix a point $x_0$. \emph{The complete shadow\/} of a point $x_0\/$ \emph{with respect to point\/} $y \in \overline{X}\setminus \{x_0\}\/$ is by definition a set \[\sha_y(x_0) := \{z \in \overline{X}|\quad \exists [yz] \quad x_0 \in [yz]\}\, .\]
Assuming of existence is necessary here only if both points $y,z\/$ are infinite: $y,z \in \partial_\infty X$.
\emph{Spherical shadow\/} of a point  $x_0\/$ \emph{of radius\/} $\rho \, >\, 0\/$ relatively point $y\in \overline{X}\/$ is intersection $\sha_y (x_0, \rho)\/$ of its complete shadow $\sha_y(x_0)\/$ with the sphere $S(x, \rho)$. In particular, when $\rho\, =\, +\infty\/$  \[\sha_y (x_0, +\infty) := \partial_\infty(\sha_y (x_0)) := \sha_y(x_0) \cap \partial_\infty X\, .\]
\end{definition}

The horosphere  $\mathcal{HS}_{y, z}\/$ where $z \in X\/$ is such a point that $\beta_y(z)\, -\, \beta_y(x_0)\, =\, \rho$, will be denoted as $\mathcal{HS}_{y, \rho}$. Here $\beta_y\/$ is arbitrary Busemann function centered in $y\in \partial_\infty X$.
Following properties of shadows are obvious.

\begin{enumerate}\label{props}
\item Sets $\sha_y(x_0)\cup \{x_0\}\/$ and $\sha_y (x_0, \rho)\/$ are closed in $\overline{X}\/$ for all $\rho \, >\, 0$;
\item $\sha_y (x_0)\, =\, \overline{\bigcup\limits_{\rho \, >\, 0}\sha_y (x_0, \rho)}\setminus \{x_0\}$;
\item If $y \in X$, then $\sha_y (x_0, \rho)\, =\, S(y, |xy|\, +\,\rho) \cap S(x, \rho)\/$ for all $\rho \, >\, 0$;
\item If $y \in \partial_\infty X$, then $\sha_y (x_0, \rho)\, =\, (\mathcal{HS}_{y, \rho}) \cap S(x, \rho)\/$  for all $\rho \, >\, 0$;
\item \label{null} If $\angle_{x_0}(y,z)\, =\, 0$, then $\sha_y(x_0)\, =\, \sha_z(x_0)$.
\end{enumerate}
The statement \ref{null} implies consequence:
\begin{Cor}
If the direction of the ray $c\,\colon[0, \infty)\, \to\, X\/$ in a point $x_0\/$ with $c(0)= x_0\/$ and $c(|xy|)\, =\, y$, has unique inverse in $x_0$, then for any two $z', z'' \in \sha_y(x_0)\/$  \[\sha_{z'}(x_0)\, =\, \sha_{z''}(x_0)\, .\]
\end{Cor}

For sufficiently small $\varepsilon \, >\, 0\/$  $\mathcal N_\varepsilon (\sha_y(x_0, \rho))\/$ will denote \emph{$\varepsilon$-neigh\-bour\-hood\/}  of spherical shadow $\sha_y(x_0, \rho)\/$
\begin{itemize}
\item in the sphere $S(y, |xy|\, +\,\rho)$, if $y \in X$, or
\item in the horosphere $\mathcal{HS}_{y, \rho}$, if $y \in \partial_\infty X$.
\end{itemize}


\begin{Th}\label{finite}
For all $x_0 \in X, y \in \overline{X}\setminus\{x_0\}, 0\, <\, \rho \, <\, +\infty\/$ and $\varepsilon \, >\, 0 \/$ there exists  $\delta \, >\, 0\/$ such that for any $x_1 \in B(x_0, \delta)$, satisfying equality $|yx_1|\, =\, |yx_0|\/$ (or $\beta_y(x_1)\, =\, \beta_y(x_0)\/$ when $y \in \partial_\infty X$), inclusion \[\sha_y(x_1, \rho) \subset \mathcal N_\varepsilon (\sha_y(x_0, \rho))\] holds.
\end{Th}

\begin{proof} 
We prove the statement in the case  $y \in X$. Situation when $y \in \partial_\infty X\/$ is similar. Assume that for any $\delta \, >\, 0\/$ there exists a point $x_\delta \in B(x_0, \delta) \cap S(y, \, |yx_{0}|)\/$ and a point $z_\delta \in S(y, |yx_0|\, +\,\rho) \setminus \mathcal N_\varepsilon (\sha_y(x_0, \rho))$, for which $x_\delta \in [yz_\delta]$. We take a sequense $\delta_n \, \to\, 0\/$ and corresponding sequences of points $x_{\delta_{n}}\/$ and $z_{\delta_{n}}$. Then $x_{\delta_{n}} \, \to\, x_0$. Since the space $X\/$ is finitely compact, one may choose converging subsequence from $z_{\delta_{n}}$. We may assume that sequence $z_{\delta_{n}}\/$ converges to a point $z \in S(y, |yx_0|\, +\,\rho)\/$ itself. Parametrizing the segment $\gamma\, =\, [yz]\/$ by arclength, we obtain 
\begin{equation}\label{sered}
|x_0 \gamma(|yx_0|)|\, \le\, |x_0 x_{\delta_{n}}| \, +\, |x_{\delta_{n}} \gamma(|yx_0|)\, \le\, \delta_n \, +\, |z_{\delta_{n}} z|.
\end{equation}
The value on the right hand of \eqref{sered} vanishes when $n \, \to\, \infty$, hence the constant on the left hand is 0. Consequently, the point  $z \in \sha_y(x_0, \rho)$, and points $z_n\/$ belong to $\mathcal N_\varepsilon (\sha_y(x_0, \rho))\/$ when  $n\/$ is sufficiently large. Contradiction. 
\end{proof}

Let a set $\mathcal V \subset \partial_\infty X\/$ and a number $K, \varepsilon \, >\, 0\/$ be given. \emph{$(y, K, \varepsilon)$-neigh\-bour\-hood\/} of $\mathcal V\/$ is by definition a set \[\mathcal N_{y,\, K,\, \varepsilon}(\mathcal V) := \{\zeta \in \partial_\infty X | \exists \xi \in \mathcal V \quad \zeta \in \mathcal U(\xi, y ,K, \varepsilon)\}\, ,\] where
\[ \mathcal U(\xi,x_0,K, \delta) := \{\eta \in \partial_\infty X |\,\, |\,c(K)d(K)|\,\, <\,\varepsilon, \quad  c=[y, \xi] ; d\, =\, [y, \eta] \}\, .\]

Following claim in fact is only reformulated statement of the theorem  \ref{finite} for a case $y \in X$.
\begin{Cor}\label{infty}
{\sloppy
For any point  $y \in X\/$ and any $(y, K, \varepsilon)$-neighbourhood $\mathcal N_{y,\, K,\, \varepsilon}(\partial_\infty(\sha_y (x_0)))\/$ of shadow at infinity $\partial_\infty (\sha_y(x_0))\/$ there exists  $\delta\, >\,0$, such that for any  $x_1\in B(x_0, \delta)$, with $|yx_1|\, =\, |yx_0|$, we have
\[\partial_\infty (\sha_y(x_1)) \subset \mathcal N_{y,\, K,\, \varepsilon} (\partial_\infty (\sha_y(x_0)))\, .\]

}
\end{Cor}

\subsection{The geometry of ideal boundary of $CAT(0)$-space}\label{Z}
\label{existence}

In this subsection we recall some well-known facts from asymptotic geometry of Hadamard spaces. We refer the reader to \cite{BH} for more detailed considerations.

First of all, we note that given a point $y_0 \in X$, an ideal point $\xi \in \partial_\infty X\/$ and a number $t \, >\, 0$, a family of sets 
\begin{equation}\label{sphebaza}
\mathcal B_{y_{0},\, \xi} := \{ \mathcal U_{\delta, t}(y_0,\, \xi)|\quad \delta,t \, >\, 0\}
\end{equation}
form a basis of neighbourhoods of point $\xi\/$ in the cone topology in $\partial_\infty X$.
Here
\[\mathcal U_{\delta, t}(y_0,\, \xi) := \{ \eta \in \partial_\infty X |\quad |c(t) d(t)|\, <\,\delta \}\, ,\]
and $c, d\,\colon [0, +\infty) \, \to\, X\/$ are rays emanating from $y_0\/$ in  $\xi, \eta \in \partial_\infty X\/$ directions correspondingly.

The space $X\/$ instead of the cone topology on ideal boundary $\partial_\infty X\/$ has one induced by so called \emph{angle metric}. For  $\xi, \eta \in \partial_\infty X\/$ angle distance is by definition equal  to  \[\angle(\xi, \eta) := \sup \{\angle_x (\xi, \eta)| \quad x \in X\}\, .\] The interior metric on  $\partial_\infty X\/$ associated with angle metric is called  \emph{Tits metric\/} and denoted as $Td$. The two metrics on $\partial_\infty X\/$ are equivalent in the sence that they induce the same topology on $\partial_\infty X$. We will write $\partial_T X\/$ for denoting the ideal boundary equiped with Tits metric.
We need following two propositions.

\begin{Prop}[\cite{BH}, Proposition 9.5]
Angle metric considered as a function $(\xi, \eta) \, \to\, \angle (\xi, \eta)\/$ is lower semicontinuous with respect to cone topology: for all
 $\varepsilon \, >\, 0\/$ there exist neighbourhoods  $\mathcal U\/$ of point  $\xi\/$ and $\mathcal V\/$ of point $\eta$, such that for all $\xi\,' \in \mathcal U\/$ è $\eta\,' \in \mathcal V\/$ inequality \[\angle (\xi\,', \eta\,') \, >\, \angle (\xi, \eta)\, -\, \varepsilon\] holds.  
\end{Prop}

As a corollary, the Tits metric is also lower semicontinuous in the cone topology. 

\begin{Prop}[\cite{BH}, Proposition 9.21]
Let $X\/$ be a proper $CAT(0)$-space and let $\xi_0\/$ and $\xi_1\/$ be distinct points of $\partial_\infty X$.
\begin{enumerate}
\item If $Td(\xi_0, \xi_1 \, >\, \pi\/$ then there is a geodesic $c\,\colon \R \, \to\, X\/$ with $c(+\infty)\, =\, \xi_0\/$ and $c(-\infty)\, =\, \xi_1$,
\item If there is no geodesic $c\,\colon \R \, \to\, X\/$ with $c(+\infty)\, =\, \xi_0\/$ and $c(-\infty)\, =\, \xi_1$, then $Td(\xi_0, \xi_1)\, =\, \angle(\xi_0, \xi_1)\/$ and there is a geodesic segment in $\partial_T X\/$ joining $\xi_0\/$ and $\xi_1$,
\item If $c\,\colon \R \, \to\, X\/$ is a geodesic, then $Td(c(-\infty), c(+\infty)) \, \ge\, \pi\/$ with equality iff $c\/$ bounds a flat half-plane,
\item If the diameter of the Tits boundary $\partial_TX\/$ is $\pi$, then every geodesic line in $X\/$ bounds a flat half-plane.
\end{enumerate}
\end{Prop}

\subsection{Points with uniqueness of inverse direction}

Directions $\xi, \eta \in \Sigma_xX\/$ are called \emph{mutually inverse\/} if $\angle_x(\xi, \eta)\, =\, \pi$. In the case of geodesically complete $CAT(0)$-space two directions $\xi, \eta\in \Sigma_xX\/$ are mutually inverse iff there exists a geodesic throw $x\/$ whose positive direction in $x\/$ is $\xi\/$ and negative one is $\eta$.

Given geodesic $a\/$ in $CAT(0)$-space $X\/$ we denote $\omega^+\,(a)\/$ a set of points $x \in a\/$ with the property that the positive direction $\xi \in \Sigma_xX:\,  [xa(+\infty)] \in \xi\/$ of $a\/$ has more than one inverse direction. Similarly $\omega^-(a)\/$ be the set of points where the negative direction $\eta \in \Sigma_xX:\,  [xa(-\infty)] \in \eta\/$ of $a\/$ has more than one inverse.

\begin{Th}\label{singpoints}
Let $a\/$ be a geodesic in a geodesically complete locally compact $CAT(0)$-space $X$.
Then sets $\omega^+(a)\/$ and $\omega^-(a)\/$ are at most countable.
\end{Th}

\begin{proof}

For $\phi \, >\, 0\/$ consider set $\Omega^+_\phi(a) \subset \omega^+(a)\/$ defined as following. $x \in \Omega^+_\phi(a)\/$ iff there exists a direction $\zeta \in \Sigma_xX\/$ inverse to direction of the ray $[xa(+\infty)]$, such that $\angle_x(\zeta, a(-\infty)) \, >\, \phi$.  We will show that intersection of $\Omega^+_\phi(a)\/$ with any segment $[xy]\subset a\/$ is finite. 

Really, assume that there exists a segment $[xy]\subset a\/$ containing infinite sequence $\{a(t_n)\}_{n=1}^{\infty}\subset \Omega^+_\phi(a)$. We may assume that $x\, =\, a(0)\/$ and $y\, =\, a(-L)$, where $L\, =\, |xy|$. For a point $a(t_n)\/$ we set $z_n \in S(x, 2L)\/$ a point such that $a(t_n) \in [xz_n]\/$ and $\angle_{a(t_n) }(a(-\infty), z_n) \, >\, \phi$. Then for $n \ne k\/$ we have $|z_n z_k| \, >\, 2L\sin\frac{\phi}{2}\/$ and the sequence $\{z_n\}_{n=1}^{\infty}\/$ does not contain any fundamental subsequence, contradicting finitely compactness of $X$.

Because of $\omega^+(a)\, =\, \bigcup\limits_{\phi \, >\, 0}\Omega^+_\phi(a)\/$ and $\Omega^+_{\phi}(a) \subset \Omega^+_{\psi}(a)\/$ for $\psi \, <\, \phi$, we have the claim for $\omega^+(a)$. Consideration of $\omega^-(a)\/$ is similar.

\end{proof}

\begin{remark} Otsu and Tanoue in \cite{OT} introduce the following notion. For $\delta \, >\, 0\/$ and $y \in X\/$ point $x\/$ is called \emph{$\delta$-branched point of $y$\/} if the diameter of the set $\{v\in \Sigma_xX| \angle_x(v_{xy}v)=\pi\}\/$ is not smaller then $\delta$. Here $v_{xy}\in \Sigma_xX\/$ is the direction of the segment $[xy]$. Such a notion has evident extension on the case $y \in \overline{X}$. In fact, every set $\Omega^+_\phi(a)\/$ is contained in a set of $\phi$-branching points of $a(+\infty)\/$ lying in $a$, which also has finite intersection with every segment $[xy]\subset a$.
\end{remark}

\begin{example}
Every point of a geodesic $a\subset X\/$ may occur to be a branching point. For example, exclude from Euclidean plane the interior domain bounded by a parabola and attach a flat half-plane on its place. We get a $CAT(0)$-space $X$. The parabola $a\/$ dividing two flat domains is its geodesic and every its point is a point of branching of geodesics in both directions. However we have $\omega^+(a)\, =\, \omega^-(a)\, =\, \emptyset$, since angle between any two branches of $a\/$ in its opposite direction starting from the same point vanishes.
\end{example}

\begin{example}
Sets $\omega^+(a)\/$ and $\omega^-(a)\/$ may occur to be dence in $a$. To see this one may take a convex continuous natural parametrised curve on which every point has a positive and negative semitangents, not opposite to each other in rational points, instead of the parabola in previous example.
\end{example}

\subsection{Existence of scissors}

We prove the existence theorem for scissors in this subsection.

{\sloppy
\begin{Th} \label{15.02.03-4}
Let  $a\,\colon (-\infty , +\infty) \, \to\, X\/$ be a geodesic of strictly rank one. Let $x_0 \in a\setminus (\omega^+(a)\cup\omega^-(a))\/$ be a point where both directions of $a\/$ have unique inverse. Then there exists a geodesic $a'\/$ with $a'(0)= x_0\/$ and $\angle_{x_{0}}(a(+\infty), a'(+\infty))\, =\, 0\/$ with following property. For every neighbourhood $\mathcal U\/$ of a triple \[(a'(+\infty), a'(-\infty), x_0) \in \partial_\infty X \times \partial_\infty X \times X\] there exists a triple $(\xi,  \eta, x) \in \mathcal U\/$ with $x \ne x_0\/$ and geodesics $b'\, =\, [a'(-\infty) \xi], c'=[\eta a'(+\infty)]\/$ and $d'=[\eta\xi]$, forming scissors $\langle a',b',c',d'; x\rangle$.
\end{Th}

}
\begin{proof}

We note that by the condition on the rank of $a$, any geodesic $a'\/$ with $a'(+\infty) \in \partial_\infty(\sha_{a(-\infty)}(x_0))\/$ or $a'(-\infty) \in \partial_\infty(\sha_{a(+\infty)}(x_0))\/$ has rank $1$. First we show that there exists scissors with lowest base $a\/$ and center $x\/$ arbitrary closed to $x_0$. In view of the first remark, the considerations will be applicable also for a geodesic $a'\/$ passing throw $x_0\/$ in the same direction.

Take points $y'\, =\, a(-\rho)\/$ and $y''\, =\, a(\rho)$, where $\rho \, >\, 0$. We have 
\[\partial_\infty(\sha_{a(-\infty)}(x_0))\, =\, \partial_\infty(\sha_{y'}(x_0))\] and
\[\partial_\infty(\sha_{a(+\infty)}(x_0))\, =\, \partial_\infty(\sha_{y''}(x_0))\, .\]

Low semicontinuous funtion 
\[Td\,\colon \partial_\infty X \times \partial_\infty X \, \to\, \R_+\cup \{+\infty\}\]
attains its minimum on the compact set \[Q=\partial_\infty(\sha_{y'}(x_0)) \times \partial_\infty(\sha_{y''}(x_0))\, ,\] and inequality 
\begin{equation}\label{mintd}
{\rm min}(Td)|_Q \, >\, \pi
\end{equation}
holds as a consequence of the condition on the rank of $a$. Moreover, there exists neighbourhoods 
\[\mathcal N\,':=\mathcal N_{y', K, \varepsilon}(\partial_\infty(\sha_{y'}(x_0)))\] 
and 
\[\mathcal N\,'':=\mathcal N_{y'', K, \varepsilon}(\partial_\infty(\sha_{y''}(x_0)))\] with some $K \, >\, \rho\/$ and $\varepsilon \, >\, 0$,
for which 
\begin{equation}\label{inftd}
\inf \{Td(\xi, \eta)|\quad (\xi, \eta) \in \mathcal N'\times \mathcal N''\} \, >\, \pi.
\end{equation}

{\sloppy
Choose $\delta_1$-neighbourhood $B(x_0, \delta_1)\/$ of a point $x_0\/$ defined by Corollary \ref{infty} relatively $\mathcal N_{y', K, \varepsilon/2}(\partial(\sha_{y')}(x_0, \rho)\/$ and $\mathcal N_{y'', K, \varepsilon/2}(\partial(\sha_{y'')}(x_0, \rho)$.

}
Let also $\mathcal N_{\varepsilon/2}(\partial(\sha_{a(-\infty)}(x_0, \rho)))\/$ and $\mathcal N_{\varepsilon/2}(\partial(\sha_{a(+\infty)}(x_0, \rho)))\/$ be $\varepsilon/2$-neighbourhoods of boundaries of shadows of point $x_0\/$ relatively points $a(-\infty)\/$ and $a(+\infty)$.

From Theorem \ref{finite} there exists $\delta_2$-neighbourhood $B(x_0, \delta_2)\/$ of point $x_0$, such that for every $x' \in B(x_0, \delta_2)\/$ inclusions 
\[\sha_{a(-\infty)}(x', \rho) \subset \mathcal N_{\varepsilon/2}(\partial(\sha_{a(-\infty)}(x_0, \rho)))\]
and
\[\sha_{a(+\infty)}(x', \rho) \subset \mathcal N_{\varepsilon/2}(\partial(\sha_{a(+\infty)}(x_0, \rho)))\]
hold.

Set $\delta_0 := \min \{\delta_1, \delta_2\}$. Then for any point $x \in \mathcal U_{\delta_{0}}(x_0)\/$ and geodesics $b\/$ and $c\/$ satisfying conditions
\begin{itemize}
\item $b(0)\, =\, c(0)\, =\, x$,
\item $b(-\infty)\, =\, a(-\infty)\/$ and
\item $c(+\infty)\, =\, a(+\infty)$
\end{itemize}
we get
\begin{equation}\label{bincl}
b(\rho) \in \mathcal N_{\varepsilon/2}(\partial(\sha_{a(-\infty)}(x_0, \rho)))
\end{equation}
and
\[c(-\rho) \in \mathcal N_{\varepsilon/2}(\partial(\sha_{a(+\infty)}(x_0, \rho)))\]

We will show that
\begin{equation}\label{b+infty}
b(+\infty) \in \mathcal N\,'
\end{equation}
and
\begin{equation}\label{c-infty}
c(-\infty) \in \mathcal N\,''.
\end{equation}

Given geodesic ray $\gamma\, =\, [y'b(+\infty)]\/$ with length parametrisation \[\gamma\,\colon [0, +\infty) \, \to\, X\] and geodesic line $a'\/$ passing throw $x_0\/$ such that $a'(+\infty) \in \partial_\infty(\sha_{y'}(x_0))\/$ we have
\[|\gamma(2\rho)a'(\rho)|\, \le\, |\gamma(2\rho)b(\rho)| \, +\,|b(\rho)a'(\rho)|\, .\]

First item has estimation
\[|\gamma(2\rho)b(\rho)|\, \le\,|\gamma(0)b(-\rho)|\, =\, |a(-\rho)b(-\rho)|\, \le\, |a(0)b(0)| \, <\, \frac\varepsilon2\, .\]

Because of \eqref{bincl} one may choose the geodesic $a'\/$ such that second item satisfies to inequality
\[|b(\rho)a'(\rho)|\, <\,\frac\varepsilon2\, .\]
Finally we get
\[|\gamma(2\rho)a'(\rho)|\, <\,\varepsilon\, ,\]
proving the inclusion \eqref{b+infty}. Inclusion \eqref{c-infty} is similar.

So because of (\ref{inftd}) there exists a geodesic $d$ in $X$ connecting points $c(-\infty)$ and $b(+\infty)$ such that we have scissors $\langle a,b,c,d; x\rangle$.

{\sloppy
Now we take a sequences $\delta_n \, \to\, 0\/$ of scales and $\langle a_n, b_n, c_n, d_n; x_n \rangle\/$ of corresponding scissors for which $|x_0x_n| \, <\, \delta_n$. Choose as $\xi,\eta \in \partial_\infty X\/$ limit points of sequences $b_n(+\infty)\/$ and $c_n(-\infty)\/$ correspondingly. Then $\xi \in \partial_\infty(\sha_{a(-\infty)}(x_0))\/$ and $\eta \in \partial_\infty(\sha_{a(+\infty)}(x_0))$. The points above may be connected by a geodesic $a'=[\eta\xi]\/$ in $X\/$ such that $a'(0)=x_0\/$ (cf. inequality \eqref{mintd}). For any $\delta \, >\, 0\/$ there exists scissors $\langle a', b', c', d';x\rangle\/$ with lowest base $a'\/$ and $|xx_0| \, <\, \delta$. Note that the point $x\/$ may always be different of $x_0\/$ and not belong to $a'$.

}
It is left to show that such scissors may be chosen with $b'(+\infty)\/$ arbitrary closed to $\xi$, and $c'(-\infty)\/$ arbitrary closed to $\eta\/$ in the sence of cone topology on $\partial_\infty X$. It may be done with applying the same method as above but with neighborhoods of points $a'(\pm \infty)\/$ instead of neighborhoods of shadows at infinity of point $x_0$. Construction of a line $a'\/$ garantees that sets
\[C_{-\infty}(B(a'(K),\varepsilon))\, =\, \bigcup\limits_{y\in B(a'(K),\varepsilon)}[a'(-\infty)y]\]
and
\[C_{+\infty}(B(a'(-K),\varepsilon))\, =\, \bigcup\limits_{z\in B(a'(-K),\varepsilon)}[za'(+\infty)]\]
has nonempty intersection for any $\varepsilon, K\, >\,0$, and moreover 
\begin{align*}
C_{-\infty}&(\mathcal N_\varepsilon (a'(K))) \cap C_{+\infty}(\mathcal N_\varepsilon (a'(-K)))\\
 \cap (X \setminus (\sha_{y'} &(x_0) \cup \sha_{y''}(x_0)\cup \{ x_0\}) \cap B(x_0, \varepsilon)) \ne \emptyset.
\end{align*}
This fact provides for necessary construction.

\end{proof}

\begin{remark}\label{aa'}
Note that geodesics $a\/$ and $a'\/$ are connected by the asymptotic chain $a\, =\, a_0, a_1, a_2\, =\, a'$, where the geodesic $a_1\/$ is  obtained as union of rays $[x_0 a(-\infty)]\/$ and $[x_0 a'(+\infty)]$.
\end{remark}

\subsection{Continuity of the displacement function}

The goal of this subsection is a theorem \ref{contin} which is a theorem of continuity for a displacement function $\delta\/$ defined by the scissors translations along given geodesic $a\/$ as a function defined on appropriate subset of $\partial_\infty X \times \partial_\infty X \times X$. First we need some estimation of distance between projections.

We say that a point $x \in X\/$ \emph{projects\/} to a point $x_0 \in a$, if $x_0\/$ is  nearest to $x\/$ point of a geodesic $a$. It is called a \emph{projection\/} of $x$. Every point $x\in X\/$ has unique projection on any given geodesic $a$.

\begin{Th}\label{vareps}
Let
$d\,\colon~(-\infty,\,+\infty)~\, \to\,~X\/$ be a strongly rank one geodesic, $x \notin d\/$ be a point in $X\/$ and $x_1\/$ --- its projection onto $d$. Then for any  $\varepsilon \, >\, 0\/$ points $\xi\, =\, d(+\infty)\/$ and $\eta\, =\, d(-\infty) \in \partial_\infty X\/$ have neighbourhoods  $\mathcal U_+\/$ and $\mathcal U_-\/$ correspondingly, such that if a geodesic $d'\/$ connects points $d'(+\infty) \in \mathcal U_+\/$ and $d'(-\infty) \in \mathcal U_-$, and $x'_1 \in d'\/$ is a projection of the point $x\/$ onto $d'$, then $|x_1 x'_1| \, <\, \varepsilon$. 
{\sloppy

}
\end{Th}

The first step of its proof is the next simple lemma.

\begin{Le}\label{tech}
Let $d\/$ be a geodesic in the $CAT(0)$-space $X\/$ and $x'\/$ be the projection of the point $x\/$ onto $d$. Then for any point $y \in d\/$ \[|\, x'y|\, \le\, \sqrt{|\, xy|^2\, -\, |\, xx'|^2}\, .\]
\end{Le}

\begin{proof} Assume that $|\, x'y| \, >\, \sqrt{|\, xy|^2\, -\, |\, xx'|^2}.\/$ Consider a comparison triangle $\overline{x} \overline{x}' \overline{y}\/$ of a triangle $xx'y$. We have by assumption $\angle_{\overline{x}'}(\overline{x},\, \overline{y}) \, <\, \pi/2$. Hence there exists a point $\overline{m} \in \overline{x}' \overline{y}$, with $|\, \overline{x}\, \overline{m}| \, <\, |\, \overline{x}\, \overline{x}'|$. Its corresponding point $m\/$ of triangle $xx'y\/$ satisfies inequality $|\, xm|\, \le\, |\, \overline{x}\, \overline{m}| \, <\, |\, xx'|$. The contradiction proves the claim. 
\end{proof}

\begin{proof}[Proof of the theorem \ref{vareps}]
Let $d\,\colon~(-\infty,\,+\infty)~\, \to\,~X\/$ be given geodesic and $\sigma \, >\, 0\/$ an arbitrary number.  We denote $HS_{\xi}(t, \rho)\/$ the intersection of a horosphere $\mathcal{HS}_{\xi, d(t)}\/$ centered at a point $\xi \in \partial_\infty X\/$ with $B(d(t), \rho)$.
{\sloppy

}

We need the fact that the cone topology of  $CAT(0)$-space $X\/$ has a base consisting of sets  \[\mathcal U(\xi,x_0,K, \delta) := \{\eta \in \partial_\infty X |\,\, |\,c(K)c_*(K)|\,\, <\,\delta \quad  c=[x_0, \xi] ; c_*\, =\, [x_0, \eta] \}\]
 
Choose points $x_{-K}\, =\, d(-K)\/$ and $x_K\, =\, d(K)$, where $K\, >\,0\/$ is a sufficiently large nubmer. In particular one may assume for $K\/$ to satisfy the following condition: if  $y\in  B(x_{-K},\, \sigma) $, then $[y d(+\infty)] \cap HS_{d(+\infty)}(K,\, \sigma/4) \ne \emptyset$, and if $z\in B(x_K, \, \sigma)$, then $[z d(-\infty)] \cap HS_{d(-\infty)}(-K,\, \sigma/4) \ne \emptyset$. This condition can be accomplished because  $d\/$ has strictly rank one.{\sloppy

}

Denote as
\[\mathcal U_{+\infty}(K, \sigma) := \left\{ \xi \in \partial_\infty X |\, \forall y \in B(x_{-K}, \sigma) | \quad [y \xi] \cap HS_{d(+\infty)}\left(K,\, \frac{\sigma}{2}\right) \ne \emptyset\right\}\]
 and 
\[\mathcal U_{-\infty}(K, \sigma) := \left\{ \eta \in \partial_\infty X |\, \forall z \in B(x_K,\, \sigma) | \quad [z \eta] \cap HS_{d(-\infty)}\left(-K,\, \frac{\sigma}{2}\right) \ne \emptyset\right\} \]
sets of those ideal points which serve as centers  of projections onto corresponding horospheres, such that a ball $B(x_{-K},\, \sigma)\/$ moves inside a ball $B(x_K,\,\sigma/2)$, and a ball $B(x_K, \sigma)\/$ moves inside a ball $B(x_{-K}, \sigma/2)$. Sets $\mathcal U_{\pm\infty}(K, \sigma)\/$ are nonempty, since $d(+\infty) \in \mathcal U_{+\infty}(K, \sigma)\/$ and $d(-\infty) \in \mathcal U_{-\infty}(K, \sigma)$. Furthermore, points $\gamma(\pm \infty)\/$ are interior points of sets  $\mathcal U_{\pm\infty}(K, \sigma)$. This is a consequence of the definition of the cone topology on $\partial_\infty X\/$ and local compactness of $X$. 

For an arbitrary point $x \notin d\/$ let $x_1\/$ be its projection onto $d$. Fix $\varepsilon \, >\, 0$. One may assume that $\varepsilon \, <\, \min \{|\, x \, x_1|, \, 1 \}$, besause it is sufficient to make all considerations for arbitrary small $\varepsilon$. Set 
\[\sigma := \frac{\varepsilon^2}{9|\, x\, x_1|}\, .\] 
For  $K\/$ as above we denote $\mathcal U_\pm\subset \partial_\infty X\/$ neighbourhoods of ideal points $d(\pm \infty)$, contained correspondingly in $\mathcal U_{\pm\infty}(K, \sigma)$.

We show that neighbourhoods  $\mathcal U_\pm\/$ are the ones satisfying the claim of the theorem: if the geodesic $d'\/$ has ends $d'(\pm \infty) \in \mathcal U_\pm\/$ and $x'_1\in d'\/$ is a projection of a point $x\/$ onto $d'$, then $|\, x_1\, x'_1| \, <\, \varepsilon$. For this we show that $d' \cap \overline{HS}_{d(+\infty)}(K,\, \sigma/2) \ne \emptyset\/$ and $d' \cap \overline{HS}_{d(-\infty)}(-K,\, \sigma/2) \ne \emptyset\/$ at first.

Indeed, let $\pi_+ \,\colon HS_{d(-\infty)}(-K,\, \sigma) \, \to\, HS_{d(+\infty)}(K, \, \sigma/2)\/$ be the projection map centered in $\gamma'(+\infty)\/$ and $\pi_- \,\colon  HS_{\gamma(+\infty)}(K,\, \sigma) \, \to\, HS_{\gamma(-\infty)}(-K, \, \sigma/2)\/$ be analogous projection centered in $d'(-\infty)$. Then by the condition on the rank of $d'$, the composition $\pi_- \circ \pi_+\, : HS_{d(-\infty)}(-K, \sigma) \, \to\,  HS_{d(-\infty)}(-K, \sigma/2)\/$ is a continuous map which is  contraction operator on $HS_{d(-\infty)}(-K, \sigma)$. Since the closure $\overline{HS}_{d(-\infty)}(-K, \, \sigma/2)\/$ is compact, the contraction operator  $\pi_- \circ \pi_+\/$ has unique fixed point in $\overline{HS}_{d(-\infty)}(-K, \, \sigma/2)$. Let such a point $z\/$ have an image $y\, =\, \pi_+(z)$. Then $y \in [z \, d'(+\infty)]\/$ and $z \in [y \, d'(- \infty)]$. The union of rays $[z \, d'(+\infty)] \ni y\/$ and $[y \, d'(- \infty)] \ni z\/$ gives us precisely a geodesic $d'$, since it does not bound any flat strip. 

It is only left to estimate the distance $|\, x'_1\, x_1|\/$ to complete the proof. Since  $d'\/$ passes throw interior points of balls $B(-K,\, \sigma)\/$ and $B(K,\, \sigma)$, there exists a point $m \in d'$, with $|\, m\, x_1| \, <\, \sigma$. Then we have:
\begin{equation}\label{first}
|\, x'_1\, x_1 |\, \le\, |\, x'_1\, m| \, +\, |m \, x_1|\, \le\, \sqrt{|\, x\, m|^2\, -\, |\, x\, x'_1|^2} \, +\, \sigma.
\end{equation}

By the triangle inequality we have for point $m\/$ 
\begin{equation} \label{m}
|\, x\, m| \, <\, |\, x \, x_1| \, +\, \sigma,
\end{equation}
and since  $x_1\/$ is the nearest to $x\/$  point of $d\/$ and the distance of $x'_1\/$ to $d\/$ is less then $\sigma$, the estimate 
\begin{equation}\label{x1}
|\, x \, x'_1| \, >\, |\, x \, x_1|\, -\, \sigma
\end{equation}
holds.

Substituting inequalities \eqref{m} and \eqref{x1} to \eqref{first}, we get
\begin{align*}
|\, x'_1 \, x_1| \, <\, \sqrt{(|\, x\, x_1| \, +\, \sigma)^2\, -\, (|\, x\, x_1|\, -\, \sigma)^2} \, +\, \sigma\, &=\, 2\sqrt{\sigma|\, x \, x_1|} \, +\, \sigma\\
= \sqrt{\sigma}(2 \sqrt{|\, x \, x_1|} \, +\, \sqrt{\sigma}) \, <\, 3 \sqrt{\sigma |\, x x_1|}\, &=\, \varepsilon.
\end{align*}

\end{proof}

Fix a strictly rank one geodesic $a\,\colon \R \, \to\, X$. Set $Z(a)\subset \partial_\infty X \times \partial_\infty X \times X\/$ to be a subset consisting of all triples $(\xi, \eta, x) \in \partial_\infty X \times \partial_\infty X \times X\/$ such that there exists scissors $\langle a,b,c,d; x\rangle\/$ with $b(+\infty)\, =\, \xi\/$ and $c(-\infty)\, =\, \eta$. 

\begin{Th}\label{contin}
The displacement function $\delta\/$ is continuous on the set $Z(a)$.
\end{Th}

\begin{proof} We use the equality \eqref{15.02.03-1}. Fix a triple $(\xi_0, \eta_0, x_0) \in Z(a)$.  
It means that  $x_0\/$ is the center of scissors $\langle a, b_0, c_0, d_0; x_0 \rangle$, where $b_0(+\infty)\, =\, \xi_0\/$ and $c_0(-\infty)\, =\, \eta_0$. 

Fix an arbitrary $\varepsilon \, >\, 0$. Then first, by continuity of  Busemann functions $\beta_{a-}\/$ and $\beta_{a\, +\,}\/$ there exists  $\sigma_1\/$ such that if for a point $x'\in X\/$ inequality $|\, x_0x' | \, <\, \sigma_1\/$ holds, then 
\begin{equation}\label{15.02.03-2}
|\, \beta_{a\, +\,}(x') \, +\, \beta_{a-}(x')\, -\, \beta_{a\, +\,}(x_0)\, -\, \beta_{a-}(x_0) | \, <\, 
\varepsilon/2.
\end{equation}

Using the theorem \ref{vareps} we choose neighbourhoods $\mathcal U_+\/$ and $\mathcal U_-\/$ of points $d_0(+\infty)=\xi_0\/$ and $d_0(-\infty)=\eta_0 \in \partial_\infty X\/$ with the condition: for any geodesic $d\,'\/$ such that $d\,'(-\infty) \in \mathcal U_-\/$ and $d\,'(+\infty) \in \mathcal U_+$,  the projection  $x_1\/$ of $x_0\/$ onto $d_0\/$ and projection $x'_1\/$ of  $x_0\/$ onto $d\,'\/$ satisfy inequality $|\, x_1x'_1 | \, <\, \varepsilon/4$. For these neighbourhoods there exists  $\sigma_2\/$ such that if the point $x'\/$ satisfy conditions\par
a) there exists a ray $[x'\xi']$, inverse to the ray $[x' a(-\infty)]\/$ with $\xi' \in \mathcal U_+$,\par
b) there exists a ray $[x'\eta']$, inverse to the ray $[x' a(+\infty)]\/$ with $\eta' \in \mathcal U_-$,\par
c) $d\,'\,\colon (-\infty, +\infty) \, \to\, X\/$ is a geodesic connecting points $\xi'\/$ and $\eta'$, and \par
d) $|\, x_0x' | \, <\, \sigma_2$,\\
then  
\[|\, \beta_{d\,'-}(x')\, -\, \beta_{d_{0}-}(x_0)| \, <\, \varepsilon/4\]
and
\[|\, \beta_{d\,'+}(x')\, -\, \beta_{d_{0}\, +\,}(x_0)| \, <\, \varepsilon/4\, .\]

{\sloppy
Here $\beta_{d'-}\/$ and $\beta_{d'\, +\,}\/$ are Busemann functions corresponding to points $d'(-\infty)\/$ and $d'(+\infty)\/$ and satisfying the condition $\beta_{d'-}(x'_1)\, =\, \beta_{d'\, +\,}(x'_1)\, =\, 0$, and $\beta_{d_{0}-}$, $\beta_{d_{0}\, +\,}\/$ --- analogous Busemann functions for the geodesic $d_0$. As a result, for $\sigma\, =\, \min \{\sigma_1 , \sigma_ 2\}\/$ and for all $x'\/$ such that $|\, x_0x' | \, <\, \sigma\/$ and conditions a), b), c) above holds, we have an estimation:
\begin{equation} \label{15.02.03-3}
|\, \beta_{d'-}(x') \, +\, \beta_{d'\, +\,}(x')\, -\, \beta_{d_{0}-}(x_0)\, -\, \beta_{d_{0}\, +\,}(x_0)|\, \le\, \varepsilon/2.
\end{equation}

}
Taking $\mathcal N\, =\, \mathcal U_+ \times \mathcal U_- \times B(x_0, \sigma)\/$ as neighbourhood of triple $(\xi_0, \eta_0, x_0)\/$ and compairing inequalities \eqref{15.02.03-2} and \eqref{15.02.03-3} with representation \eqref{15.02.03-1}, we get a condition of the continuity for the function $\delta\/$ as function of triple $(\xi, \eta, x)\/$ in the point  $(\xi_0, \eta_0, x_0)\in Z(a)$. \end{proof}

\begin{Th}
Let $x_0\in a\/$ be an arbitrary point of a geodesic  $a\/$ of strictly rank one. Then   $\delta(\xi, \eta, x)\, \to\, 0\/$ when $(\xi, \eta, x)\, =\, (b(+\infty), c(-\infty), x)\/$ tends to $(a(+\infty), a(-\infty), x_0)\/$ in the sence of topology in the set $Z(a)\/$ inherited from $\partial_\infty X\times \partial_\infty X\times X$.
\end{Th}

\begin{proof} Let  $A\/$ be a set of lines $a'\/$ which passes throw $x_0\/$ such that $a'(+\infty) \in \sha_{a(-\infty)} (x_0, +\infty)\/$ and $a'(-\infty) \in \sha_{a(+\infty)} (x_0, +\infty)$.
For $a' \in A\/$ we define  "\emph{closed\/}" scissors $\langle a', a', a', a'; x_0 \rangle\/$ as a system of four items of a line $a'\/$ and a center $x_0$. They also have a transformation map $T\, =\, {\rm id}_{a'}$, for which $\delta T\, =\, 0\/$ and a function  $\delta\/$ remains continuous when it is defined on  $Z(a) \cup \{(a'(+\infty), a'(-\infty), x)|~x\in a',\, a'\in A\}\/$ by equality $\delta(a'(+\infty), a'(-\infty), x)\, =\, 0\/$ for $x \in a'$. 
\end{proof}

Following evident corollary will be the key point in remainder of the proof.

{\sloppy
\begin{Cor} \label{small}
Let $a\/$ be the geodesic of strictly rank one. Then there exists a point $x_0\in a$, a geodesic $a'$, connected with $a\/$ by the asymptotic chain, and $\Delta \, >\, 0$, such that for every $\varepsilon \in (0, \Delta)\/$ there exist scissors $\langle a', b, c, d; x \rangle$, with displacement of scissors translation $T\/$ equal to $\delta T\, =\, \varepsilon$.
\end{Cor}

}

\subsection{Recovery of the metric on geodesic of strictly rank one}

The goal of this paragraph is to show that the metric of an arbitrary geodesic of strictly rank one is restorable from the diagonal tube $V$. 

\begin{Le}
Let $x \notin a\/$ and $\xi, \eta \in \partial_\infty X$. Then the relation  $V\/$ allows to detect whether the triple  $(\xi, \eta, x)\/$ belongs to a set $Z(a)$.
\end{Le}

\begin{proof} $(\xi, \eta, x) \in Z(a)\/$ if and only if there exists  $r$-sequences $\{x_z\}_{z\in \Z}\subset a$,  $\{ u_z \}_{z \in \Z}$, \, $\{ v_z \}_{z \in \Z}\/$ and $\{ w_z \}_{z \in \Z}$, for which $u_0\, =\, v_0\, =\, x$,  and such that their limiting ponts in $\partial_\infty X\/$ satisfy equalities $u_{-\infty}\, =\, x_{-\infty}$, \, $v_{+\infty}\, =\, x_{+\infty}$, \, $w_{-\infty}\, =\, v_{-\infty}=\eta\/$ and $w_{+\infty}\, =\, u_{+\infty}=\xi$. Given any four  $r$-sequences $\{ x_z \}_{z \in \Z}$, \,  $\{ u_z \}_{z \in \Z}$, \, $\{ v_z \}_{z \in \Z}\/$ and $\{ w_z \}_{z \in \Z}$, relation  $V\/$ detects, whether marked equalities hold or not. Also, $V\/$ lets to reveal such four  $r$-sequences if they exist. \end{proof}

\begin{Le}
Given scissors $\langle a,\, b,\, c,\, d;\, x \rangle$, relation $V\/$ defines an image $T(m)\in a\/$ of any point $m \in a\/$ in scissors translation $T$.
\end{Le}

\begin{proof} We use the fact that all horospheres are defined by $V\/$ and serves as level sets of Busemann functions. Hence  $V\/$ defines a point $m_1\, =\, R_{ac}(m)\/$ as $\mathcal{HS}_{a(+\infty),\, m} \cap c$, a point $m_2\, =\, R_{cd}(m_1)\/$ as $\mathcal{HS}_{d(-\infty),\, m_1} \cap d$, a point  $m_3\, =\, R_{db}(m_2)\/$ as $\mathcal{HS}_{d(+\infty),\, m_2} \cap b\/$ and finally, a point $T(m)\, =\, m_4\, =\, R_{ba}(m_3)\/$ as $\mathcal{HS}_{a(-\infty),\, m_3} \cap a$. 
\end{proof}

\begin{Le}\label{sdvig}
Let $\langle a,\, b,\, c,\, d;\, x \rangle\/$ be scissors with displacement of translation $\delta T$. Then  $V\/$ allows to display the value $\delta T$.
\end{Le}

\begin{proof}  We have:
\[
 \delta T\, =\, \frac 1n d(x_0, T^n(x_0))\, =\, \lim\limits_{n \, \to\, +\infty} \frac 1n d(x_0, T^n(x_0))\, =\, \lim\limits_{n \, \to\, +\infty} \frac 1n [d(x_0, T^n(x_0))],
\]
where $[\, t]\/$ is integral part of the number $t$. Since every item of the sequence $a_n\, =\, [d(\, x_0, T^n(x_0)\, )]\/$ is definable by  $V$, hence the value $\delta T\/$ is. 
\end{proof}

We are ready now to prove the following theorem.

\begin{Th}\label{hyperbolic}
Let metric space $(X,d)\/$ and trial metric $d'\/$ on $X\/$ be as in Theorem \ref{tube}.
Assume that $a \,\colon (-\infty , \, +\,\infty) \, \to\, X\/$ is a geodesic of strictly rank one in metric $d$. Then $a\/$ is of strictly rank one in metric $d'\/$ as well, and for all $t_1, t_2 \in \R\/$ 
\begin{equation}\label{dd2}
d'(c(t_1), c(t_2))\, =\, d(c(t_1), c(t_2)).
\end{equation}
\end{Th}

\begin{proof} 
We may think the assertion on the rank of $a\/$ as already proved  in Proposition \ref{rank}, since the property of geodesic to have rank one is accomplished when any its $r$-sequence is of rank one.

According Corollary \ref{small} there exists $x_0 \in a$, geodesic $a'\/$ and $\Delta\, >\,0\/$ such that for any $q \in \N$, satisfying $1/q \, <\, \Delta$, there exist scissors $\langle a',\, b,\, c,\, d;\, x \rangle\/$ with displacement $\delta T\, =\, 1/q$.  By lemma \ref{sdvig} the relation $V\/$ lets to detect, whether the value of displacement equal $1/q$. Take scissors with this value of displacement, which is independent on the choise of metric $d\/$ or $d'$. We have \[d(a'(t), T^p(a'(t)))\, =\, d'(a'(t), T^p(a'(t)))\, =\, p/q\]
for all $t \in \R$.
Described procedure defines all points of type $a'(t)\/$ with $t \in \Q\/$ in $a'\/$ from relation $V$, and we have $d(a'(0), a'(t))\, =\, d'(a'(0), a'(t))\/$ for them. When $t \in \R \setminus \Q\/$ we may use the fact that $V\/$ defines the incidence relation for $a'\/$ and its order of points.

It remains to apply remark \ref{aa'} and Lemma \ref{asymchain} for completing the proof.
\end{proof}

Combining Theorem \ref{hyperbolic} with Theorem \ref{parabolic}, one gets the whole proof of the Theorem \ref{tube}. Consequently  Theorem \ref{main} is also proved.

\end{document}